\documentclass{birkjour}
\usepackage{amssymb}
\usepackage{amsthm}

\usepackage[colorlinks]{hyperref}
\usepackage{color,graphicx,shortvrb}

\usepackage[active]{srcltx} %SRC Specials for DVI search

\usepackage{enumerate}
\usepackage{amssymb}
\usepackage{tikz-cd}

\usepackage{mathtools}

 \newtheorem{theorem}{Theorem}[section]
 \newtheorem{corollary}{Corollary}[section]
 \newtheorem{lemma}{Lemma}[section]
 \newtheorem{proposition}{Proposition}[section]
 \theoremstyle{definition}
 
 \theoremstyle{remark}
 \newtheorem{remark}{Remark}[section]

 \numberwithin{equation}{section}

\begin{document}

%-------------------------------------------------------------------------
% editorial commands: to be inserted by the editorial office
%
%\firstpage{1} \volume{228} \Copyrightyear{2004} \DOI{003-0001}
%
%
%\seriesextra{Just an add-on}
%\seriesextraline{This is the Concrete Title of this Book\br H.E. R and S.T.C. W, Eds.}
%
% for journals:
%
%\firstpage{1}
%\issuenumber{1}
%\Volumeandyear{1 (2004)}
%\Copyrightyear{2004}
%\DOI{003-xxxx-y}
%\Signet
%\commby{inhouse}
%\submitted{March 14, 2003}
%\received{March 16, 2000}
%\revised{June 1, 2000}
%\accepted{July 22, 2000}
%
%
%
%---------------------------------------------------------------------------
%Insert here the title, affiliations and abstract:
%

\title[$M$-ideals, yet again: the case of real JB$^*$-triples]{$M$-ideals, yet again: the case of real JB$^*$-triples}

\author[D.P. Blecher]{David P. Blecher}

\address{Department of Mathematics, University of Houston, Houston, TX 77204-3008, USA.}
\email{dblecher@math.uh.edu}
%\thanks{The research was supported by balabala.}

\author[M. Neal]{Matthew Neal}
\address{%
Math and Computer Sci. Department, Denison University, Granville, OH 43023,  USA.}
\email{nealm@denison.edu}
%\thanks{The research was supported by balabala.}

\author[A.M. Peralta]{Antonio~M.~Peralta}
\address{Instituto de Matem{\'a}ticas de la Universidad de Granada (IMAG), Departamento de An{\'a}lisis Matem{\'a}tico, Facultad de
	Ciencias, Universidad de Granada, 18071 Granada, Spain.}
\email{aperalta@ugr.es}
%\thanks{The research was supported by balabala.}

\author[S. Su]{Shanshan Su}
\address{%
School of Mathematics, East China University of Science and Technology, Shanghai, 200237 China. \\
(Current address) Departamento de An{\'a}lisis Matem{\'a}tico, Facultad de
	Ciencias, Universidad de Granada, 18071 Granada, Spain.}
\email{lat875rina@gmail.com}
%\thanks{The research was supported by China Scholarship Council.}

%----------classification, keywords, date


\date{today}
%----------additions
%\dedicatory{To my boss}
%%% ----------------------------------------------------------------------

\begin{abstract}
We prove that a subspace of a real JBW$^*$-triple is an $M$-summand if and only if it is a weak$^*$-closed triple ideal. As a consequence, $M$-ideals of real JB$^*$-triples correspond to norm-closed triple ideals. As in the setting of complex JB$^*$-triples, a geometric property is characterized in purely algebraic terms. This is a newfangled treatment of the classical notion of $M$-ideal in the real setting by a fully new approach due to the unfeasibility of the known arguments in the setting of  complex C$^*$-algebras and JB$^*$-triples. The results in this note also provide a full characterization of all $M$-ideals in real C$^*$-algebras, real JB$^*$-algebras and real TROs.      
\end{abstract}

%%% ----------------------------------------------------------------------
\maketitle
%%% ----------------------------------------------------------------------
%\tableofcontents
\section{Introduction}\label{sec1}

It seems quite a long time since E.M. Alfsen and E.G. Effros introduced in 1972 the %purely geometric 
notions of $M$- and $L$-summands and $M$-ideals of an abstract Banach space \cite{AlfEff72}. However, we must admit, almost forty years later, that these objects have not been fully described in some non-trivial (but natural) structures.  Let us begin by fixing the central notions in this note. \smallskip

Let $V$ be a Banach space. A projection $P$ on $V$ is called an \emph{$L$-projection} (respectively, an \emph{$M$-projection}) if  $$\|  x \| = \| P(x) \| + \| (Id-P)(x) \|, \quad \hbox{ for all } x \in V$$ (respectively, $\|  v \| = max \{ \| Pv \| , \| (Id-P)(v) \|  \}, \quad \forall v \in V)$.  The image of an $L$-projection (respectively, an $M$-projection) on $V$ is called an \emph{$L$-summand} (respectively, an \emph{$M$-summand}).  A closed subspace $M$ of a Banach space $V$ is said to be an \emph{$M$-$ideal$} if $M^{\circ}:=\{\varphi\in V^* : \varphi|M \equiv 0\}$ is an $L$-$summand$ of $V^*$, that is, there is an $L$-projection $Q$ on $V^*$ whose image is $M^{\circ}$. In such a case $N := (Id-Q) (V^*)$ is a closed subspace of $V^*$ with $V^* = M^{\circ} \mathop{\oplus}\limits^{\ell_{1}} N$. Moreover, there exists a closed subspace $J$ of $V^{**}$ such that $V^{**} = J \mathop{\oplus}\limits^{\ell_{\infty}} M^{**}$, the direct sum of $J$ and $M^{**}$ with the $\ell_{\infty}$-norm, where we identify $M^{**}$ with the weak$^*$-closure of $M$ in $V^{**}$. The undoubted attraction of these notions in Banach space theory, and close disciplines such as approximation theory, has motivated an intensive study of these notions in several classes of Banach spaces. The book by P.~Harmand, D.~Werner, and W.~Werner \cite{HarmandWerners1993}  is widely mentioned as the reference source on the general M-ideal theory and their applications, despite the appearance of subsequent works on this topic.\smallskip

Though Alfsen and Effros' original aim was to employ the theory of $M$-ideals in the setting of C$^*$-algebras, ordered Banach spaces and $L_1$-preduals, they defined the central notions only in terms of the geometry of the Banach space without any algebraic or order theoretic structure. However, in the case of C$^*$-algebras, JB$^*$-algebras and JB$^*$-triples,  $M$-ideals can be characterized in purely algebraic terms. More concretely, R.R. Smith and J.D. Ward showed in \cite{SmidWard78} that $M$-ideals in a Banach algebra with identity are subalgebras, which are in fact ideals if the algebra is commutative, and the $M$-ideals of a C$^*$-algebra are precisely the closed two-sided ideals (see also the historical comments in \cite[\S V.7]{HarmandWerners1993}).  It is due to Alfsen and Effros that if $A_{sa}$ is the self-adjoint part of a C$^*$-algebra $A$, then the $M$-ideals in $A_{sa}$ are just the subspaces of the form $I\cap  A_{sa},$  with $I$ being a uniformly closed two-sided ideal in $A$ (see \cite[Proposiiton 6.18]{AlfEff72}). A result by R. Payá, F.J. Pérez and Á. Rodríguez-Palacios asserts that the weak$^*$-closed ideals in a JBW$^*$-algebra correspond to its $M$-summands, and consequently, the closed ideals in a JB$^*$-algebra are in one-to-one correspondence with its $M$-ideals (see \cite{PaPeRo1982}). In the wider setting of JB$^*$-triples (see section~\ref{pre} for detail), G. Horn showed that every weak$^*$-closed triple ideal of a JBW$^*$-triple is an $M$-summand \cite[Theorem 4.2]{Horn1987structure}, while J.T. Barton and R. Timoney concluded in \cite[Theorem 3.2]{BartonTimoney1986} that the closed triple ideals in a JB$^*$-triple $E$ are precisely the $M$-summands of the Banach space $E$. It is the interplay between geometric and algebraic structures increases the interest on these notions.\smallskip

Attention to real versions of C$^*$-algebras and JB$^*$-triples has increased in recent times. Their importance has grown  specially in the study of metrics invariants and properties which are stable under surjective isometries between certain subsets of these structures, as in Tingley's problem or in the Mazur-Ulam property (see, for example, \cite{Pe2018survey, MoriOza2019, BeCuFerPe2020, KaPe21}). Surprisingly, the algebraic characterization of $M$-ideals in real C$^*$-algebras and real JB$^*$-triples as closed ideals is still missing. Further results towards the development of the theory of one-sided $M$-ideals in real structures have been developed by S. Sharma \cite[\S 5]{Sharma14}. Recall that a \emph{complete left $M$-projection} on a \emph{real operator space} (i.e. a closed subspace of $B(H)$, for some real Hilbert space $H$), $X$, is a projection $P$ on $X$ such that the map $\nu_{_P}^{c} :  X \to  C_2 (X)$, $x \mapsto \left[\begin{matrix} P(x) \\
	x - P(x)
\end{matrix} \right]$ is a complete isometry. A subspace $M$ of a real operator space $X$ is a \emph{right $M$-ideal} if $J^{\circ\circ}$ is the range of a complete left $M$-projection on $X^{**}$. A result by Sharma proves that the right $M$-ideals in a real C$^*$-algebra $A$ are precisely the closed right ideals in $A$ (see  \cite[Corollary 5.10]{Sharma14} and \cite{BleEffrZari2002} for the result in the setting of C$^*$-algebras).\smallskip

The explanation for the lack of a proof in the real setting can be easily understood by taking a simple look at the tools employed in the complex setting. Let us focus, for example, on the widest possible setting of complex JB$^*$-triples. In such a case, for every $x$ in a JB$^*$-triple $E$, the mapping $e^{i t L(x,x)}$ (where $L(x,x)$ is the bounded linear operator on $E$ given by $L(x,x) (a) = \{x,x,a\}$) is a surjective linear isometry on $E$, and hence $L(x,x)$ leaves invariant each $M$-ideal in $E$ (cf. \cite[Proof of Theorem 3.2]{BartonTimoney1986} or the arguments in \cite[Corollary 4.2]{PaPeRo1982} and Proposition~\ref{prop0} here). However, in the real setting it is more difficult to find logarithms of surjective real linear isometries, our available tools reduce to consider exponentials of the linear operators of the form $\delta (x,y) (a) = \{x,y,a\} - \{y,x,a\}$ (see Corollary~\ref{c image under an inner derivation}). There is another handicap in the real setting. Namely, Kaup's Banach--Stone theorem, which asserts that a linear bijection between complex JB$^*$-triples is a triple isomorphism if and only if it is an isometry (see  \cite[Proposition 5.5]{Ka83}), does not necessarily hold for real JB$^*$-triples (cf. \cite{ChuDaRuVen,Da,FerMarPe2004}). \smallskip

That is, the essential arguments in the setting of complex JB$^*$-triples are simply unachievable for real JB$^*$-triples. Quoting K. Goodearl,  \emph{``The change of coefficient field from $\mathbb{C}$ to $\mathbb{R}$ is more than just a cosmetic change.''} \cite{Goodearl82}. It becomes necessary to develop a complete new strategy to characterize $M$-ideals in real JB$^*$-triples. The characterization of $M$-ideals for real JB$^*$-triples is a key result to complete our knowledge on the geometric properties of these objects, and at the same time open the path to establish new decompositions of real JBW$^*$-triples as the sum of their type I and continuous parts.\smallskip

This note contains complete new arguments to prove that $M$-ideals in real C$^*$-algebras, real JB$^*$-algebras (see Corollary~\ref{c M-ideals of real JB-algebras}) and real JB$^*$-triples (cf. Theorem~\ref{thm1}) are two-sided ideals and triple ideals, respectively. We actually show in Theorems~\ref{t M.-ideals in real JBWstar algebras} and \ref{t M-summands in real JBW*-triples are ideals} that every $M$-summand in a real JBW$^*$-algebra (respectively, in a real JBW$^*$-triple) is a weak$^*$-closed ideal (respectively, triple ideal). To avoid the inherent difficulties appearing in the setting of real JB$^*$-triples, arguments rely on an appropriate application of the algebraic characterization of the facial structure of the closed unit ball of a real JB$^*$-triple and the natural connections with $M$-decompositions (see Proposition~\ref{prop1} and Theorem~\ref{t complete tripotents and faces}). %No need to say that the new arguments are also valid in the complex case.
 Our efforts to avoid these obstacles in the real setting have resulted in more-geometric proofs, which in turn could lead to more general applications.

\section{Preliminaries}\label{pre}

According to \cite{Ka83}, a complex Banach space $\mathcal{E}$ is called a (complex) \emph{JB$^*$-triple} when it can be equipped with a continuous triple product $\{\cdot, \cdot, \cdot  \}$ from $\mathcal{E} \times \mathcal{E} \times \mathcal{E}$ into $\mathcal{E}$, which is conjugate linear in the middle variable, linear and symmetric in the outer variables and satisfies the following axioms:\begin{enumerate}[$(i)$]
\item For every $x,y,z,a,b \in \mathcal{E}$, we have 
$$\{a,b, \{x,y,z \} \}= \{ \{a,b,x\},y,z \} - \{x, \{b,a,y\},z \} +\{x,y, \{a,b,z \} \} ;$$ \hfill (\emph{Jordan identity})
\item For each $a\in \mathcal{E}$, the linear operator $x \mapsto L(a,a) (x):=\{a,a,x\}$ is a hermitian operator on $ \mathcal{E}$ with non-negative spectrum;
\item $\| \{a,a,a \} \| = \|a\| ^{3}$ for all $a \in \mathcal{E}$. \hfill (\emph{Gelfand-Naimark axiom})	
\end{enumerate}

For all $x,y \in \mathcal{E}$, we shall denote by $L(x,y)$ and $Q(x,y)$ the linear and conjugate linear mappings on $\mathcal{E}$ given by 
$$ L(x,y)(a) = \{x,y,a \}, \quad Q(x,y)(a) =\{x,a,y \} \quad (\forall a \in \mathcal{E}),$$  respectively. We shall write $Q(a) $ for $Q(a,a)$.\smallskip

All structures in this note are assumed to be non-trivial. The closed unit ball of a Banach space $X$ will be denoted by $\mathcal{B}_{X}$. \smallskip

The theory of JB$^*$-triples has been developed during the last five decades. Since every $C^*$-algebra is a JB$^*$-triple with triple product \begin{equation}\label{eq triple product Cstar}  \{a,b,c \} = \frac{1}{2} (ab^{*}c + c b^{*}a), 
\end{equation} many of the basic properties and results in C$^*$-algebra theory have been pursued in the strictly wider setting of JB$^*$-triples. In some cases, the latter setting is more general and complicated to handle. The class of JB$^*$-triples also includes all JB$^{*}$-algebras with the triple product \begin{equation}
\label{eq triple product JBstar}  \{a,b,c \} = (a \circ b^{*}) \circ c + (c \circ b^*) \circ a - (a \circ c) \circ b^*.
\end{equation} For the basic background on JB$^*$-algebras the reader is referred to the monographs \cite{HOS,AlfShu2003,CabRod2014}. Let us recall that every C$^*$-algebra can be naturally regarded as a JB$^*$-algebra with respect to the natural Jordan product defined by $a \circ b := \frac{1}{2} (ab +ba)$.\smallskip

A \emph{JBW$^*$-triple} is a JB$^*$-triple which is also a dual Banach space. In a clear analogy with Sakai's theorem in the case of von Neumann algebras (cf. \cite[Theorem 1.7.8]{Sak71}), every JBW$^*$-triple admits  a unique isometric predual and its triple product is separately weak$^*$-continuous (see \cite{BartonTimoney1986}).  The bidual, $\mathcal{E}^{**},$ of each JB$^{*}$-triple $\mathcal{E}$ is a JBW$^*$-triple \cite{Di86}. \smallskip

Y. Friedman and B. Russo proved in \cite{FriRu86} that the norm and the triple product of each JB$^*$-triple $\mathcal{E}$ satisfy the following inequality: \begin{equation}\label{eq cont triple products} \|\{x,y,z\}\| \leq \|x\| \ \|y\| \ \|z\|, \ \ (x,y,z\in \mathcal{E}).
\end{equation} 

The result known as \emph{Kaup's Banach-Stone theorem} states that a linear bijection between JB$^*$-triples is an isometry if and only if it is a triple isomorphism  (cf. \cite[Proposition 5.5]{Ka83}).\smallskip

A closed JB$^*$-subtriple $\mathcal{I}$ of a JB$^*$-triple $\mathcal{E}$ (i.e. $\mathcal{I}$ is a norm-closed subspace of $E$ and $\{\mathcal{I},\mathcal{I},\mathcal{I}\}\subseteq \mathcal{I}$) is said to be an \emph{ideal}  or a \emph{triple ideal} of $\mathcal{E}$ if $\{ \mathcal{E},\mathcal{E},\mathcal{I}\} + \{ \mathcal{E},\mathcal{I},\mathcal{E}\}\subseteq \mathcal{I}$. For example, the triple ideals of a C$^*$-algebra $\mathcal{A}$ are precisely the (two-sided) ideals of $\mathcal{A}$ (see \cite[\S 4]{Horn1987structure}, \cite[\S 3]{BartonTimoney1986} and \cite[Proposition 1.3]{BuChu92} for characterizations and additional details). Along this note, an element $p$ in a C$^*$-algebra or in a JB$^*$-algebra is called a projection if $p =p^2 = p^*$. \smallskip

Real operator algebras, real C$^*$-algebras and real von Neumann algebras were studied studies in the decades of 1980's and 1990's. K. Goodearl emphasizes in \cite{Goodearl82} that \emph{``The change of coefficient field from $\mathbb{C}$ to $\mathbb{R}$ is more than just a cosmetic change.''}  A real C$^*$-algebra is a real Banach $^*$-algebra $\mathcal{A}$ satisfying  $\|a a^*\| = \|a \|^2$ and $ \mathbf{1} + aa^* $ is invertible in the unitization of $\mathcal{A}$, for all $a\in \mathcal{A}$ (in the complex case, the second axiom can be relaxed, see \cite{Goodearl82,Li2003}). It is known that $\mathcal{A}$ is a real C$^*$-algebra if, and only if, any of the following equivalent statements holds: 
\begin{enumerate}[$(\checkmark)$]
\item The complexification of $\mathcal{A}$ is a C$^*$-algebra under an appropriate norm.
\item There exists a C$^*$-algebra $A$ and a conjugate-linear $^*$-automorphism $\tau$ on $A$ such that $\mathcal{A} = A^{\tau} =\{x\in A : \tau(x) =x\}$.
\item $\mathcal{A}$ is isometrically $^*$-isomorphic to a uniformly closed $^*$-subalgebra of $B(H)$ for some real Hilbert space $H$. 
\end{enumerate}

\emph{Real JB$^*$-triples} were originally defined as those real norm-closed JB$^*$-subtriples of (complex) JB$^*$-triples (see \cite{IsidroKaupPalacios}). Therefore the class of real JB$^*$-triples contains all complex JB$^*$-triples, and all real and complex C$^*$-algebras. As in the case of real C$^*$-algebras, real JB$^*$-triples can be obtained as \emph{real forms} of JB$^*$-triples. More concretely, given a real JB$^*$-triple $E$, there exists a unique (complex) JB$^*$-triple structure on its algebraic complexification $\mathcal{E}= E \oplus i E,$ and a conjugation (i.e. a conjugate linear isometry of period 2) $\tau$ on $\mathcal{E}$ such that $$E = \mathcal{E}^{\tau} = \{ z\in \mathcal{E} : \tau (z) = z\},$$ (see \cite{IsidroKaupPalacios}). Let us observe that by  Kaup's Banach-Stone theorem, the mapping $\tau$ preserves triple products.  In particular, the real Banach space $B(H_1,H_2)$ of all bounded real linear operators between two real, complex, or quaternionic Hilbert spaces is also a real JB$^*$-triple with the triple product given in \eqref{eq triple product Cstar}. Real JB$^*$-algebras (i.e. real norm-closed JB$^*$-subalgebras of JB$^*$-algebras) are also included in the class of JB$^*$-triples when they are equipped with the triple product in \eqref{eq triple product JBstar} (cf. \cite{Alvermann,Pe2003ax} and \cite[\S 4]{RoPa2010}). Clearly, every real JB$^*$-triples also satisfy the Jordan identity and the inequality 
\begin{equation}\label{eq cont triple products real} \|\{x,y,z\}\| \leq \|x\| \ \|y\| \ \|z\|, \ \ (x,y,z\in \mathcal{E}),
\end{equation}  also holds for the norm and the triple product of any real JB$^*$-triple.\smallskip 

A real JBW$^*$-triple is a real JB$^*$-triple which is additionally a dual Banach space. It is known that, as in the complex case, each real JBW$^*$-triple admits a unique isometric predual and its triple product is separately weak$^*$-continuous \cite{MartinezPeralta2000}.\smallskip

One handicap of real JB$^*$-triples is that Kaup's Banach-Stone theorem is no longer, in general, valid on these structures. That is, there are surjective linear isometries between real JB$^*$-triples which are not triple isomorphisms. The best we can, in general, say is that any such a mapping preserves the symmetrized triple product given by $\langle x,y,z \rangle = \frac13 \left(\{x,y,z\} + \{y,z,x\} + \{z,x,y\}\right)$ (cf. \cite[Theorem 4.8]{IsidroKaupPalacios}). Surjective linear isometries between real JB$^*$-triples have been intensively studied and under certain additional hypotheses on the involved real JB$^*$-triples (for example, when both of them are real JB$^*$-algebras) they are triple isomorphisms (cf. \cite{ChuDaRuVen,Da,FerMarPe2004}).\smallskip

An element $e$ in a real or complex JB$^*$-triple ${E}$ is called a \emph{tripotent} if $\{e,e,e\}= e$. The symbol $\mathcal{U} ({E})$ will stand for the set of all tripotents in ${E}$. Every tripotent $e\in E$ gives rise to the following decompositions of $E$ $$ E = E_{0} (e) \oplus E_{1} (e) \oplus
E_{2} (e) = E^{0} (e) \oplus  E^{1} (e) \oplus E^{-1} (e),$$ where
$E_{k} (e) := \{ x\in E : L(e,e)x = \frac{k}{2} x \}$ is a
subtriple of $E$ and $E^{k} (e) := \{ x\in E : Q(e)
(x) := \{ e,x,e \}= k x \}$ is a real Banach subspace of $E$ (compare \cite[Theorem 3.13]{Loos2}). The natural projections of $E$ onto $E_{k} (e)$ and $E^{k} (e)$ will be denoted by $P_{k} (e)$ and
$P^{k} (e)$, respectively. The first decomposition is known as the
\emph{Peirce decomposition} of $E$ relative to $e$, and the corresponding summands are called the \emph{Peirce}-$k$ subspaces. The canonical projection of $E$ onto ${E}_{k}(e)$ is called the Peirce $k$-projection of ${E}$ for $k=0,1,2$. Triple products among elements in the Peirce subspaces obey  certain patterns known as \emph{Peirce rules}:
$$\begin{aligned}
\{ {E_{i}(e)},{E_{j}(e)},{E_{k}(e)} \} & \subseteq E_{i-j+k}(e),\ \hbox{
if  } i,j,k\in\{0,1,2\} \\ 
\hbox{ and }   \{ {E_{i}(e)},{E_{j}(e)}, &{E_{k}(e)} \} =\{0\}, \hbox{ for } i-j+k\neq 0,1,2. \\
\{ {E_0 (e)},{E_2 (e)},{E}\} &= \{ {E_2 (e)},{E_0 (e)},{E }\} =\{0\}
\end{aligned}$$ The subspaces in the second decomposition also satisfy the following rules:
\begin{displaymath}
	\begin{array}{c}
		E_{2} (e) = E^{1} (e) \oplus E^{-1} (e), \ \ \ \ \ E_{1} (e)
		\oplus E_{0} (e) = E^{0} (e) \ \ \ \ \ \ \ \ \ \ \ \vspace{0.3cm} \nonumber\\
		\{ {E^{i} (e)},{E^{j} (e)},{E^{k} (e)}\} \subseteq E^{i j k} (e),
		\hbox{ whenever } i j k \ne 0. \ \ \ \ \ \ \ \ \ \ \ \nonumber
\end{array}\end{displaymath}

If $\mathcal{E}$ is a (complex) JB$^*$-triple and $e\in \mathcal{U} (\mathcal{E})$, the Peirce-2 subspace $\mathcal{E}_{2}(e)$ is a JB$^*$-algebra with identity $e$, Jordan product $x\circ_e y :=\{x,e,y\}$ and involution $x^{*_e} = \{e,x,e\}$.  Furthermore, $\mathcal{E}^{-1} (e) = i \mathcal{E}^{1} (e)$. In case that $E$ is a real JB$^*$-triple,  $E_2 (e)$ is a unital real JB$^*$-algebra with respect to the given operations, and its self-adjoint part, $E^{1} (e)$, is a JB-algebra.\smallskip

All Peirce projections are contractive (i.e. $\|P_k(e)\|\leq 1$ for all $k=0,1,2$, see \cite[Corollary 1.2$(a)$]{FriedmanRusso1985}). The algebraic expressions defining Peirce projections read as follows:  $P_{2}(e)= Q(e)^2$, $P_{1}(e) = 2(L(e,e) - Q(e)^{2})$, $P_{0}(e) = Id - P_{1}(e)-P_{2}(e) = Id- 2L(e,e) + Q(e)^2$.\smallskip

A tripotent $e$ in a real or complex JB$^*$-triple $E$ is called \emph{complete} if ${E}_{0}(e) =\{0\}$. We shall denote by $\mathcal{U}_{c} ({E})$ the set of all complete tripotents in ${E}$. Another interesting connection between the algebraic and geometric properties of JB$^*$-triples asserts that the extreme points of the closed unit ball of a real or complex JB$^*$-triple $E$ ($\partial_e \left(\mathcal{B}_{E}\right)$ in short) are precisely the complete tripotents in $E$, that is, $\partial_e \left(\mathcal{B}_{E}\right) = \mathcal{U}_{c} ({E})$ (see \cite[Lemma 4.1]{BraKaUp} and \cite[Lemma 3.3]{IsidroKaupPalacios}). In particular, by the Krein–Milman theorem, every real JBW$^*$-triple contains a broad set of complete tripotents. \smallskip

Tripotents $e,v$ in a real or complex JB$^*$-triple $E$ are called \emph{compatible} if the Peirce projections associated with $e$ and $v$ commute, that is, $P_j (e) P_k (v) = P_k(v) P_j (e)$ for all $j,k\in \{0,1,2\}$. It is known that if $e\in E_{j}(v)$ for some $j\in\{0,1,2\}$, the tripotents $e,v$ are compatible (cf. \cite[$(1.10)$]{Horn1987structure}). When $e$ and $v$ are compatible we have a simplified joint Peirce decomposition associated to the pair $e,v$ given by $\displaystyle E = \bigoplus_{j,k= 0,1,2} \big(E_j(e) \cap E_k (v)\big)$, where the natural projection of $E$ onto $E_{j,k} = E_j(e) \cap E_k (v)$ is $P_j (e) P_k (v) = P_k(v) P_j (e)$.

\subsection{Orthogonality and $M$-orthogonality}

Let $x$ and $y$ be elements in a real or complex JB$^*$-triple $E.$,We say that $x$ and $y$ are \emph{orthogonal} ($x\perp y$ in short) if $L(x,y)=0$ (see \cite[Lemma 1]{BurFerGarMarPe} for additional characterizations). In particular, two tripotents $e$ and $f$ in $E$ are orthogonal if and only if $e\in E_{0} (f)$. By employing the notion of orthogonality we can define a partial ordering on $\mathcal{U} (E)$ defined by $e\leq v$ if and only if $v-e \in \mathcal{U} (E)$ with $v-e \perp e$ (equivalently, $e$ is a projection in the unital JB$^*$-algebra $E_2(v)$ \cite[Lemma 3.3]{FriedmanRusso1985}). Observe that two projections $p,q$ in a real JB$^*$-algebra $A$ are orthogonal if and only if $p\circ q =0$.\smallskip

It is known that any two orthogonal elements $x,y$ in a real JB$^*$-triple $E$ are geometrically $M$-orthogonal (i.e. $\|x \pm y \| = \max\{\|x\|,\|y\|\}$) (see \cite[Lemma 1.3$(a)$]{FriedmanRusso1985}). Let us note that the reciprocal implication is not, in general, true.\smallskip

As in the complex setting, a closed real JB$^*$-subtriple $I$ of a real JB$^*$-triple ${E}$ is said to be an \emph{ideal} or a \emph{triple ideal} of $E$ if $\{ E,E,I\} + \{ E,I,E\}\subseteq I$. We also recall that a closed subspace $I$ of a real JB$^*$-algebra $A$ is a (Jordan) ideal of $A$ if $I\circ A \subseteq I$. It is known that triple ideals of $A$ are precisely the Jordan ideals of $A$, and they are all self-adjoint (see, for example, the discussion in \cite[page 8]{GarPe2021}).\smallskip

The above geometric implications are enough to conclude that every triple ideal in a real JB$^*$-triple is an $M$-ideal. 

\begin{remark}\label{r wstar triple ideals are M-summands} Let $I$ be a triple ideal of a real JB$^*$-triple $E$. It follows from the separate weak$^*$-continuity of the triple of $E^{**}$ \cite{MartinezPeralta2000} that $I^{**}\cong \overline{I}^{w^*}$ is a triple ideal of $E^{**}$. If $\mathcal{E}$ denotes the complex JB$^*$-triple obtained by complexifying of $E$, it is known that $\mathcal{E}^{**}$ identifies with the complexification of $E^{**}$ (cf. \cite[Proof of Lemma 4.2]{IsidroKaupPalacios}, more details on the dual spaces will be given in the next subsection). It is easy to check that the complexification of $I$, $\mathcal{I} = I\oplus i I$, is a triple ideal of $\mathcal{E}$, whose second dual, $\mathcal{I}^{**}$, identified with its weak$^*$-closure in $\mathcal{E}^{**},$ is a weak$^*$-closed ideal of $\mathcal{E}^{**}.$ It follows from \cite[Theorem 4.2$(4)$ and Lemma 4.4]{Horn1987structure} that $\mathcal{I}^{**}$ is an $M$-summand in $\mathcal{E}^{**}$, and hence $I^{**}$ is an $M$-summand in $E^{**}$. This assures that $I$ is an $M$-ideal of $E$. We have shown that
\begin{equation}\label{eq w*-c-ideals are M-summands}\begin{aligned}
\hbox{Every}& \hbox{ weak$^*$-closed triple ideal of a real JBW$^*$-triple} \\
&\hbox{is an $M$-summand.}	
\end{aligned}
\end{equation}	
\begin{equation}\label{eq ideals are M-ideals}\hbox{Every triple ideal of a real JB$^*$-triple is an $M$-ideal.}
	\end{equation}
\end{remark}

The next technical result is also related to the notion of orthogonality and partial ordering, when the latter is understood in a wider setting. 

\begin{lemma}{\rm(\cite[Lemma 1.6]{FriedmanRusso1985} and \cite[Lemma 5.1]{PoloFranciscoPeralta2018})}\label{lemma3}  
Let ${E}$ be a real {\rm(}or complex{\rm)} JB$^*$-triple. Suppose $e\in E$ is a tripotent and $x\in E$ is norm-one element such that $P_{2}(e)(x) = e$. Then we have 
	$$ x= P_{2}(e)(x) + P_{0}(e)(x) = e + P_{0}(e)(x).$$
\end{lemma}

\subsection{Orthogonality in the predual space of a real JBW$^*$-triple}

According to the standard terminology, a functional $\phi$ in the predual of a (real or complex) JBW$^*$-algebra $M$ is said to be \emph{faithful} if for each $a\geq 0$ in $M$, $\phi (a) =0$ implies $a=0$.\smallskip

Given a non-zero functional $\varphi$ in the predual of a JBW$^*$-triple $\mathcal{W}$ and a tripotent $e\in \mathcal{W}$, it is known that $\varphi = \varphi P_2(e)$ if and only if $\|\varphi\| = \| \varphi|_{\mathcal{W}_2 (e)}\|$ \cite[Proposition 1]{FriedmanRusso1985}). Furthermore, by \cite[Proposition\ 2]{FriedmanRusso1985}, there exists a unique tripotent $e= s(\varphi) \in \mathcal{W}$ (called the \emph{support tripotent} of
$\varphi$, and denoted by $s(\varphi)$) satisfying $\varphi =
\varphi P_{2} (e),$ and $\varphi|_{\mathcal{W}_{2}(e)}$ is a faithful normal positive functional on the (complex) JBW$^*$-algebra $\mathcal{W}_{2} (e)$.\smallskip

Given $\phi$ and $\psi$ in $\mathcal{W}_*$, we say that $\phi$ and $\psi$ are \emph{orthogonal} ($\phi \perp \psi$) if their support tripotents are orthogonal. Elements $x,y$ in a Banach space $X$ are said to be \emph{$L$-orthogonal} ($x \bot_{L} y$ in short) if $\| x \pm y \| = \| x \| + \| y \|$. Given two subsets $\mathcal{S}_1$ and $\mathcal{S}_2$ we write $\mathcal{S}_1 \perp_{L} \mathcal{S}_2$ if and only if $x\perp_{L} y$ for all $x\in \mathcal{S}_1$ and $y\in \mathcal{S}_2$.\smallskip

In the predual of a complex JBW$^*$-triple $\mathcal{W}$ (algebraic) orthogonality is equivalent to geometric $L$-orthogonality, that is, for any $\varphi,\phi\in \mathcal{W}_*$ we have 	$\phi \bot_{L} \psi$ if and only if $ s(\phi) \bot s(\psi)$ (see \cite[Lemma\  2.3]{FriRu87} and \cite[Theorem 5.4 and Lemma 5.5]{Morthogonal}).\smallskip

The setting of real JBW$^*$-triples deserves an independent treatment. First, suppose that $W$ is a real JBW$^*$-triple regarded as a real form of a complex JBW$^*$-triple $\mathcal{W}$, that is, $W = \mathcal{W}^{\tau}$, where $\tau$ is a period-2 conjugate linear isometry on $\mathcal{W}$. Since $\tau$ must be weak$^*$-continuous (see \cite[Proposition 2.3]{MartinezPeralta2000}), the assignment $\varphi\mapsto \tau^{\sharp}$, with $\tau^{\sharp} (\varphi) (x) := \overline{\varphi (\tau(x))}$ defines another period-2 conjugate linear isometry on $\mathcal{W}_{*}$. For each $\phi \in \left( \mathcal{W}_{*} \right)^{\tau^{\sharp}}$ we have $\phi (W)\subseteq \mathbb{R}.$ It is easy to see that the mapping $\phi \in \left( \mathcal{W}_{*} \right)^{\tau^{\sharp}} \to \phi|_{W} \in W_*$ is a surjective linear isometry between the corresponding spaces\label{dual space of a real form} (cf. \cite{MartinezPeralta2000,EdwardsRuttimann2001}).\smallskip

Given $\phi$ in $W_*$ and a tripotent $e\in W$ such that $\phi (e) = 1 = \|\phi\| =1$, then $\phi = \phi P^{1} (e)$ (see \cite[Lemma 2.7]{PeSta2001}). Furthermore, since $\phi|_{W_2(e)}$ is a normal state of the real JBW$^*$-algebra $W_2(e)$ we also have \begin{equation}\label{eq positive normal states and hermitian values}\phi (x) = \phi P_2(e) (x) = \phi P^1 (e) (x) = \phi \{e,x,e\} = \phi \left(P_2(e)(x)^{*_e}\right),
\end{equation}	for all $x\in W$.  Actually, if we regard $\phi$ as a $\tau^{\sharp}$-symmetric functional in $\mathcal{W}_*$, it follows from the uniqueness of the support tripotent of $\phi$ in $\mathcal{W}$ that $\tau (s(\phi)) = s(\phi)\in W = \mathcal{W}^{\tau}$. Therefore, $s(\phi)$ will be called the \emph{support tripotent} of $\phi$ in $W$. Actually $s(\phi)$ is the unique tripotent in $W$ such that $\phi = \phi P^{1} (e)$ and $\phi|_{W^{1}(e)}$ is a faithful normal state on $W^1 (e)$. \smallskip

Since the second dual of a real JB$^*$-triple is a real JBW$^*$-triple \cite[Lemma 4.2]{IsidroKaupPalacios}, all the previous considerations hold for the dual space of every real JB$^*$-triple.\smallskip

The next lemma is stated in \cite{ApazoglouPeralta2013} for functionals in the dual space of a real JB$^*$-triple but the arguments remain valid for predual spaces of real JBW$^*$-triples.  

\begin{lemma}\cite[Lemma 3.6]{ApazoglouPeralta2013} \label{lemma1}
Let $\phi, \psi$ be normal functionals in  the predual of a real JBW$^*$-triple $W$. Then $\phi$ and $\psi$ are $L$-orthogonal if and only if they are orthogonal in $W_*$, that is,
	$\phi \bot_{L} \psi$ if and only if $ s(\phi) \bot s(\psi).$  
\end{lemma}

\subsection{Facial structure of the closed unit ball}

Let $F$ and $C$ be convex subsets of a real or complex Banach space $X$ with $F\subseteq C$. We say that $F$ is a face of $C$ if for any $x_{1}, x_{2} \in C$ and $t \in (0,1)$, the condition $(t x_{1} + (1-t) x_{2}) \in F$ implies that $x_{1}, x_{2} \in F$. An interesting case appears when $C$ is the closed unit ball of $X$. According to this notation, a point $x_0\in \mathcal{B}_{X}$ is an extreme point of $\mathcal{B}_{X}$ if and only if the set $\{x_0\}$ is a face of $\mathcal{B}_{X}$. The facial structure of closed unit ball of certain Banach spaces has been intensively studied and completely described in the following cases: 
\begin{enumerate}[$(\checkmark)$] \item Weak$^*$-closed faces of the closed unit ball of a JBW$^*$-triple and norm-closed faces of the closed unit ball of its (isometrically unique) predual space (C.M. Edwards and G.T. Rüttimann \cite{EdRutt88}).
\item Norm-closed faces of the closed unit ball of a C$^*$-algebra and weak$^*$-closed faces of the closed unit ball of its dual space (C.A. Akemann and G.K. Pedersen \cite{AkPed92}). They also rediscovered the result of C.M. Edwards and G.T. Rüttimann in the case of von Neumann algebras, which are particular examples of JBW$^*$-triples. 
\item Norm-closed faces of the closed unit ball of a JB$^*$-triple (C.M. Edwards, F.J. Fernández-Polo, C. Hoskin, and A.M. Peralta \cite{EdFerHosPe2010}), and weak$^*$-closed faces of the closed unit ball of its dual space (F.J. Fernández-Polo and A.M. Peralta \cite{FerPe10}).
\item Weak$^*$-closed faces of the closed unit ball of a real JBW$^*$-triple and norm-closed faces of the closed unit ball of its (isometrically unique) predual space (C.M. Edwards and G.T. Rüttimann \cite{EdwardsRuttimann2001}).
\item Weak$^*$-closed faces of the closed unit ball of a JBW-algebra and norm-closed faces of the closed unit ball of its (isometrically unique) predual space (M. Neal \cite{Neal2000}). 
\item Norm-closed faces of the closed unit ball of a real JB$^*$-triple and weak$^*$-closed faces of the closed unit ball of its dual space (M. Cueto and A.M. Peralta \cite{CuPe2019}). 
\end{enumerate}

The usefulness of the previous studies relies on the characterization of norm-closed faces of the closed unit ball (a property purely geometric) in algebraic terms thanks to the tripotents elements in the involved real or complex JB$^*$-triple. For our purposes in this note we recall the concrete statement in the case of norm-closed faces of the closed unit ball of the predual of a real JBW$^*$-triple. \smallskip

Let us first fix some of the standard notation (see \cite{EdRutt88, EdwardsRuttimann2001, EdFerHosPe2010, CuPe2019}). Let $V$ be a Banach space. For each couple of subsets $A \subseteq \mathcal{B}_{V}$ and $B\subseteq \mathcal{B}_{V^{*}}$, we set 
$$ A^{\prime} = \{ \varphi \in \mathcal{B}_{V^{*}}: \varphi (x) = 1 \quad \forall x \in A \} , B_{\prime} = \{ x\in \mathcal{B}_{V}: \varphi (x) = 1  \quad \forall \varphi \in B \}.$$
Obviously, $A^{\prime}$ is a norm-closed face of $\mathcal{B}_{V^*}$ and $B_{\prime}$ is a weak$^{*}$-closed face of $\mathcal{B}_{V}$. We say that $F$ is a \emph{norm-semi-exposed face} of $\mathcal{B}_{V}$ (respectively, $G$ is a \emph{weak$^*$-semi-exposed face} of $\mathcal{B}_{V^*}$) if $F=(F^{\prime})_{\prime}$ (respectively, $G= (G_{\prime})^{\prime}$). It is known that the mappings $F \mapsto F^{\prime}$
and $G \mapsto G_{\prime}$ are anti-order isomorphisms between the
complete lattices $\mathcal{S}_n(\mathcal{B}_{V})$, of norm-semi-exposed faces
of $\mathcal{B}_V,$ and $\mathcal{S}_{w^*}( \mathcal{B}_{V^*}),$ of weak$^*$-semi-exposed
faces of $\mathcal{B}_{V^*}$, and are inverses of each other.\smallskip

When $W$ is a real JBW$^*$-triple with predual $W_*$, the norm-closed faces in $\mathcal{B}_{W_*}$ are described in the following result due to Edwards and Rüttimann.

\begin{theorem}\label{t ER nclosed faces predual ball}\cite[Theorem 3.7]{EdwardsRuttimann2001} Let $W$ be a real JBW$^*$-triple. Then every norm-closed proper face of $\mathcal{B}_{W_{*}}$ is norm-semi-exposed, and the mapping
	%\begin{equation}\label{eq order isomorphism norm-closed faces predual complex} 
$u \mapsto \{u\}_{\prime}$
%\end{equation}
 is an order isomorphism from the set ${\mathcal
		U}(W)\backslash \{0\}$ of all non-zero tripotents in $W$ onto the complete lattice $\mathcal F_{n}(\mathcal{B}_{W_{*}})$ of norm-closed proper faces of $\mathcal{B}_{W_{*}}$.
\end{theorem} 

We state next some results related to the facial structure of the closed unit ball of a real JBW$^*$-triple. 

\begin{lemma}\label{lemma4}
Let $W$ be a real JB$^*$-triple. Suppose $x \in W$ with $\| x\| = 1$ and $v$ is a tripotent in $W$ such that $P^{1} (v) (x) = v$. Then $x = v + P_{0}(v)(x)$. Consequently, if $W$ is a real JBW$^*$-triple and $\{ v \}_{_{'}} \subseteq \{ x \}_{_{'}}$, then we also have $x = v + P_{0}(v)(x)$.
\end{lemma}

\begin{proof} We can clearly assume that $v$ is non-zero. Suppose first that $v= P^{1} (v)(x)$. Let $\mathcal{W}$ denote the (complex) JBW$^*$-triple obtained by complexifying $W$. Clearly $\mathcal{W}_2 (v)$ is a JBW$^*$-algebra with unit $v$ and contains the element $x_2 := P_2 (v) (x)$, whose self-adjoint part is $\frac12 (x_2 + x_2^{*_{v}})=P^{1} (v)(x) =v.$ Clearly, $x_2$ has norm-one. Since, by the Shirshov-Cohn theorem \cite[Theorem 2.4.14]{HOS}, the JB$^*$-subalgebra of $\mathcal{W}_2 (v)$ generated by $x_2$ and $v,$ as the unit element of $\mathcal{W}_2 (v)$, is a JB$^*$-subalgebra of a C$^*$-algebra $A$ (we can further assume that $v$ is the unit of $A$), we can suppose that $x_2$ is a norm-one element in the unital C$^*$-algebra $A$ whose self-adjoint part is $\mathbf{1}$. The product of $A$ is denoted by juxtaposition and the involution by $^*$. If we write $x_2 = \mathbf{1} + i k$, with $k$ self-adjoint, we have $$ 1 = \|x_2 \|^2 = \| x_2 x_2^* \| = \| \mathbf{1} + k^2 \|, $$ which implies that $k=0$, and hence $P_2 (v) (x) = x_2=\mathbf{1} = v $. Lemma~\ref{lemma3} implies that $x = v + P_0(v) (x)$. \smallskip	
	
For the second statement, we observe that $W^{1}(v)$ is a real JBW-algebra (just apply that, by the separate weak$^*$-continuity of the triple product of $W$, the projections $P_j(v)$ and $P^{k} (v)$ are weak$^*$-continuous). Observe that $\{v\}_{_{'}}$ is nothing but the set of all normal states of this JBW-algebra. By assumptions and \cite[Lemma 2.7]{PeSta2001}, for each $\varphi\in \{v\}_{_{'}}$ we have $$\varphi (v-P^{1}(v)(x)) = \varphi P^{1}(v) (v-x) = \varphi (v-x) =0.$$ Since the normal states of a JBW-algebra separate the points \cite[\S 4]{HOS}, we derive that $v = P^{1} (v) (x).$  It follows from the first part of the proof that $x = v + P_{0}(v)(x)$.
\end{proof}

We shall also need for later purposes the next version of \cite[Lemma 5.5]{Morthogonal} in the setting of real JBW$^*$-triples. 

\begin{lemma}\label{lemma2}
Let $v$ and $w$ be tripotents in a real JBW$^*$-triple $W$. Then $v \perp w$ if and only if $\{ v \}_{_{'}}  \perp_{L}  \{ w \}_{_{'}}  $ in $W_*$.
\end{lemma}

\begin{proof} Assume first that $v\perp w$, then it clearly follows that for each $\phi \in \{ v \}_{_{'}}$ and $\psi \in \{ w \}_{_{'}}$, by the definition of support tripotents we have $s(\phi) \leq v$ and $s(\psi) \leq w$. Since $v$ and $w$ are orthogonal, it follows from the properties of the partial order that $s(\phi) \perp s(\psi),$ and thus $\phi \perp_{L} \psi$ by Lemma~\ref{lemma1}.\smallskip
	 
Suppose now that $\{ v \}_{_{'}} = \{ v \}_{_{'}}^{W_*}  \perp_{L}  \{ w \}_{_{'}} = \{ w \}_{_{'}}^{W_*} $ in $W_*$. As we commented before, if $W$ is a complex JBW$^*$-triple the conclusion was obtained in \cite[Lemma 5.5]{Morthogonal}. Let $\mathcal{W}$ denote the complexification of $W$, and let $\tau$ be a conjugate linear isometry of period-2 on  $\mathcal{W}$ satisfying $W = \mathcal{W}^{\tau}$ and $W_* = \left(\mathcal{W}_* \right)^{\tau^{\sharp}}$ (see the explanation in page~\pageref{dual space of a real form}).\smallskip

If $v\not\perp w$ in $W$, then $v\not\perp w$ in $\mathcal{W}$ (actually both statements are equivalent). By applying \cite[Lemma 5.5]{Morthogonal} we deduce that $\{ v \}_{_{'}}^{\mathcal{W}_*}  \not\perp_{L}  \{ w \}_{_{'}}^{\mathcal{W}_*}$. Therefore, there exist $\varphi_1 \in \{ v \}_{_{'}}^{\mathcal{W}_*} $ and $\varphi_2\in \{ w \}_{_{'}}^{\mathcal{W}_*}$ such that $\|\varphi_1 + \sigma \varphi_2 \|<2$ for some $\sigma \in \{\pm 1\}$. Since $\tau (v) =v$ and $\tau (w) =w$, it is not hard to see that $\tau^{\sharp} (\varphi_1) \in \{ v \}_{_{'}}^{\mathcal{W}_*}$ and  $\tau^{\sharp} (\varphi_2) \in \{ w \}_{_{'}}^{\mathcal{W}_*},$ and thus $\frac{\varphi_1 + \tau^{\sharp} (\varphi_1) }{2}  \in \{ v \}_{_{'}}^{{W}_*}$ and $\frac{\varphi_2 + \tau^{\sharp} (\varphi_2) }{2}  \in \{ w \}_{_{'}}^{{W}_*}$. We can also check that 
$$\left\| \frac{\varphi_1 + \tau^{\sharp} (\varphi_1) }{2} + \sigma \frac{\varphi_2 + \tau^{\sharp} (\varphi_2) }{2} \right\| \leq \left\| \frac{\varphi_1 + \sigma \varphi_2 }{2} \right\| + \left\| \tau^{\sharp} \left(\frac{\varphi_1 + \sigma \varphi_2 }{2}\right) \right\|<2,$$ which contradicts that $\{ v \}_{_{'}}^{W_*}  \perp_{L}  \{ w \}_{_{'}}^{W_*} $.
\end{proof}

\section{Facial structure in connection with $M$-summands of real JBW$^*$-triples}\label{kernel}

We begin with a general technical result which is probably part of the folklore on $M$-ideals in Banach spaces, but we do not know an explicit proof in the setting of real Banach spaces.\smallskip

Henceforth, given a subspace $M$ of a Banach space $V$, we shall write $M^{\circ}$ for its \emph{polar} in $V^*$, that is, $M^{\circ} =\{\varphi \in V^* : \varphi|_{M} \equiv 0\}.$ If $Y$ is a subspace of $V^*$, the symbol $Y_{\circ}$ will stand for its \emph{prepolar} in $V$ given by $Y_{\circ} =\{ x \in V : \varphi (x) =  0, \ \forall \varphi \in Y\}.$     

\begin{proposition}\label{prop0} 
Let $V$ and $X$ be real {\rm(}or complex{\rm)} Banach spaces.  Suppose $\tilde{M}$ and $\tilde{N}$ are closed subspaces of $X$ such that $X = \tilde{M}\oplus^{\ell_{1}} \tilde{N}$. Let $R: X\to X$ be a bounded linear operator satisfying that the mapping $e^{t R} : X \rightarrow X$ is a surjective linear isometry for all $t \in \mathbb{R}$. Then $R \left(\tilde{M}\right) \subseteq \tilde{M}$, $R \left(\tilde{N}\right) \subseteq \tilde{N}$, $R^{*} \left(\tilde{M}^{\circ}\right) \subseteq \tilde{M}^{\circ}$ and $R^{*} \left(\tilde{N}^{\circ}\right) \subseteq \tilde{N}^{\circ}$. Consequently, if $D: V \rightarrow V$ is a bounded linear operator such that $e^{tD} : V \rightarrow V$ is a surjective linear isometry for all $t \in \mathbb{R}$ and $M \subseteq V$ is an $M$-ideal, we have $D(M) \subseteq M$.
\end{proposition}

\begin{proof} Let $Q$ denote the $L$-projection of $X$ onto $\tilde{M}$. It follows from this fact and the assumptions that $e^{t R} Q e^{-t R}$ is also an $L$-projection on $X$. A classical result by Cunningham (see \cite[Lemma 2.2]{Cunningham1960}) proves that all $L$-projections commute, and thus we have
	$$ [e^{tR} Q e^{-tR}, Q] = 0.$$
	By considering the coefficient of $t$ in the power series of the left hand side term, we obtain $R Q- Q R Q - Q R Q+ Q R =[R Q- Q R, Q] = 0$, equivalently,
	$$ R Q + Q R = 2Q R Q.$$
	Multiplying by $Q$ on the right and left hand side of this equation, respectively, we obtain 
	$$ R Q = Q R Q = Q R,$$
	which, in particular, implies that $R(\tilde{M}) \subseteq \tilde{M},$ $R(\tilde{N}) \subseteq \tilde{N}$, $R^{*} \left(\tilde{M}^{\circ}\right) \subseteq \tilde{M}^{\circ}$ and $R^{*} \left(\tilde{N}^{\circ}\right) \subseteq \tilde{N}^{\circ}$.\smallskip
	
For the final statement suppose that $M$ is an $M$-ideal of $V$. Then there exists an $L$-projection $Q$ on $V^*$ such that $Q(V^*) = M^{\circ}$ and  $$V^{*} = M^{\circ}\oplus^{\ell_1} (Id-Q) (V^*).$$ The transpose of $D$, $D^*$, is a linear and continuous map on $V^*$ with $e^{tD^*} = \left(e^{tD}\right)^*$ being a surjective linear isometry on $V^*$. The conclusion in the first statement with $R = D^*$ and $X = V^*$ assures that  $D^*(M^{\circ}) \subseteq M^{\circ}$, $D^{**}(M^{\circ\circ}) \subseteq M^{\circ\circ}$, and thus $D(M) \subseteq M$. 
\end{proof}

\begin{remark}\label{r M-summands are weak*-closed} Suppose $W$ is a real JBW$^*$-triple and $M$ and $N$ are two $M$-summands of $W$ with $W= M \oplus^{\ell_\infty} N$. Then $M$ and $N$ are weak$^*$-closed. Indeed, let $Q$ denote the projection of $W$ onto $M$. Then, by the hypotheses, the mapping $ x\mapsto T(x) :=Q(x) - (Id-Q) (x)$ ($x\in W$) is a surjective real linear isometry on $W$. Proposition 2.3 in \cite{MartinezPeralta2000} implies that $T$ is weak$^*$-continuous, and hence the statement is clear.  
\end{remark}

For any two elements $a,b$ in a real JB$^*$-triple $E$, the symbol $\delta(a,b)$ will stand for the (bounded) linear operator on $E$ defined by $\delta (a,b) = L(a,b)-L(b,a)$.

\begin{corollary}\label{c image under an inner derivation} Let ${E}$ be a real or complex JB$^*$-triple. Suppose $M$ is an $M$-ideal of ${E}$. Then \begin{equation}\label{eq Mideals stable under delta} \delta(a,b) (M) \subseteq M \hbox{ for every } a,b\in {E}.
	\end{equation} In particular, if $M$ is an $M$-summand of a real JBW$^*$-triple $W$, we have $\delta (a,b) (M) \subseteq M$ for all $a,b\in W$.  
\end{corollary}

\begin{proof} It is known, and easily deduced from the Jordan identity and the continuity of the triple product \eqref{eq cont triple products}, that, for any $a,b \in \mathcal{E},$ the mapping $\delta(a,b) = L(a,b) - L(b,a)$ is a continuous triple derivation. This assures that $e^{t \delta}$ is a surjective linear isometry for all $t\in \mathbb{R}$ (cf. \cite[Corollary 4.8]{IsidroKaupPalacios}). Proposition~\ref{prop0} gives the desired statement. 
\end{proof}

The geometric core of our arguments is presented in the next result.

\begin{proposition}\label{prop1}
Let $e$ be a tripotent in a real JBW$^*$-triple $W$. Suppose additionally that $M$ and $N$ are two $M$-summands of $W$ with $W= M \oplus^{\ell_\infty} N$. Then the sets 
$$ F_e^{M} =  M_{\circ} \cap \{ e \}_{_{'}}  \hbox{ and }  F_e^{N} =  N_{\circ} \cap  \{ e \}_{_{'}} $$
are two norm-closed faces of $\mathcal{B} _{W_*}$ satisfying:\begin{enumerate}[$(a)$]
\item $F_e^{M} \bot_{L} F_e^{N} $.
\item  $conv(F_e^{M} \cup F_e^{N}) = \{ e \}_{_{'}} $. 
\item There exist two unique tripotents $v $ and $w $ in $W$ satisfying $v\perp w$,  $F_e^{N} =  N_{\circ} \cap  \{ e \}_{_{'}}= \{ v \}_{_{'}}$ and $F_e^{M} =  M_{\circ} \cap \{ e \}_{_{'}} =\{ w \}_{_{'}}$.
\end{enumerate}

\end{proposition}

\begin{proof} Under the hypotheses of the proposition $M$ and $N$ are weak$^*$-closed and $W_* = M_{\circ} \oplus^{\ell_1} N_{\circ}$ (cf. Remark~\ref{r M-summands are weak*-closed}). We can clearly assume that all the involved summands are non-zero. Let $Q$ denote the natural $L$-projection of $W_*$ onto $M_{\circ}$.\smallskip
	
We shall first show that $F_e^{M}$ and $F_e^{N}$ are two norm-closed faces of $\mathcal{B} _{W_*}$. We shall only prove the desired statement for $F_e^{M}$, the case of $F_e^{N}$ follows by similar arguments. Take $\phi_{1},\phi_{2} \in \mathcal{B}_{W_*}$ and $t \in (0,1)$ such that $t\phi_{1} + (1-t)\phi_{2} = \varphi \in F_e^{M}$. Obviously, $\phi_{1}, \phi_{2} \in \{ e \}_{_{'}}$ since $\{ e \}_{_{'}}$ is a face of $\mathcal{B}_{W_*}$. By applying the projection $Q$ we get  
$$ t\phi_{1} + (1-t)\phi_{2}  =\varphi=  Q(t\phi_{1} + (1-t)\phi_{2}) = t Q(\phi_{1}) +(1-t)Q(\phi_{2}).$$ Evaluating the previous identity at the point $e$ we deduce that
$$ t\phi_{1}(e) + (1-t)\phi_{2}(e) = 1 = tQ(\phi_{1})(e) +(1-t)Q(\phi_{2})(e) .$$
Hence, $Q(\phi_{j})(e) = \| Q(\phi_{j}) \|= 1,$ for all $j=1,2$. Note that 
$$ 1 = \| \phi_{j} \| = \|Q(\phi_{j})\| +\| (Id-Q)(\phi_{j}) \| = 1+ \| (Id-Q)(\phi_{j}) \| \  (\forall j=1,2).$$
Thus $\phi_{1},\phi_{2} \in F_e^{M}$. $F_e^{M}$ is norm-closed because it is the intersection of two norm-closed sets. \smallskip

Take now $\psi \in \{ e \}_{_{'}}$. Observe that $\psi = Q(\psi) +(Id-Q) (\psi),$ where $Q(\psi) \in M_{\circ}$ and $(Id-Q)(\psi) \in N_{\circ}$. To avoid trivialities we can assume that both functionals, $Q(\psi)$ and $(Id-Q) (\psi)$ are non-zero. Since $Q$ is an $L$-projection, we have $$1= \| \psi \| = \| Q(\psi) \| + \| (Id-Q) (\psi) \|.$$ Set $\psi_{1} = {\| Q(\psi) \|}^{-1} Q(\psi) $ and $\psi_{2} = {\| (Id-Q)(\psi) \|}^{-1} (Id-Q)(\psi)$, then, 
$$\psi = \| Q(\psi) \| \psi_{1} + \| (Id-Q) (\psi) \| \psi_{2}.$$
Therefore, $1= \psi (e) = \| Q(\psi) \| \psi_{1}(e) + \| (Id-Q) (\psi) \| \psi_{2}(e)$ with $1= \| Q(\psi) \| + \| (Id-Q) (\psi) \|$, which implies that $\psi_{1}(e) = \psi_{2}(e) = 1$. Thus, $\psi_{1} \in  M_{\circ} \cap \{ e \}_{_{'}} = F_e^{M}$ and $\psi_{2} \in  N_{\circ} \cap  \{ e \}_{_{'}} = F_e^{N},$ and we have further shown that $conv(F_e^{M} \cup F_e^{N}) = \{ e \}_{_{'}} $.\smallskip

Statement $(a)$ follows straightforwardly from the facts $F_e^{M} \subseteq M_{\circ},$ $F_e^{N} \subseteq N_{\circ},$ and $W_* = M_{\circ} \oplus^{\ell_1} N_{\circ}$. \smallskip

Finally, if we apply Theorem~\ref{t ER nclosed faces predual ball} (cf. \cite[Theorem 3.7]{EdwardsRuttimann2001}) to the norm-closed faces $F_e^M$ and $F_e^{N}$, we deduce the existence of two tripotents $v$ and $w$ in $W$ satisfying $F_e^M = \{v\}_{_{'}}$ and $F_e^N = \{w\}_{_{'}}$. Since by $(a)$, $F_e^M = \{v\}_{_{'}} \perp_{L} F_e^N = \{w\}_{_{'}},$ Lemma~\ref{lemma2} implies that $v\perp w$.
\end{proof}

Let $v$ be a non-zero tripotent in a real JBW$^*$-triple $W$. Since $F= \mathbb{R} v$ is a JBW$^*$-subtriple of $W$ and we can easily find a norm-one functional $\phi\in F_*$ with $\phi (v) =1$, \cite[Corollary]{Bu2001} assures the existence of a norm-preserving extension of $\phi$ to a functional in $W_*$. Well, by employing this property, or by Theorem~\ref{t ER nclosed faces predual ball} we have $\{v\}_{_{'}} \neq \emptyset.$\smallskip

Recall that the complete tripotents in a real JB$^*$-triple $E$ coincide with the  extreme points of its closed unit ball $\mathcal{B}_{E}$, that is, $\mathcal{U}_c (E) = \partial_{e} (\mathcal{B}_{E})$ (see \cite[Lemma 3.3]{IsidroKaupPalacios}).

\begin{theorem}\label{t complete tripotents and faces} Let $W$ be a real JBW$^*$-triple. Suppose additionally that $M$ and $N$ are two non-zero $M$-summands of $W$ with $W= M \oplus^{\ell_\infty} N$, and let $P$ denote the natural projection of $W$ onto $M$. Then 
for each complete tripotent {\rm(}equivalently, each extreme point of the closed unit ball of $W${\rm)} $e$ in $W$, the elements $P(e)$ and $(Id-P)(e)$ are two orthogonal non-zero tripotents in $W$. Moreover, $ P(e) =v \in  \partial_e\left( \mathcal{B}_{M} \right)$ and $(Id-P) (e) = w \in \partial_e \left( \mathcal{B}_{N} \right).$ 
\end{theorem}

\begin{proof} Observe that $e= P(e) + (Id-P)(e)$ with $P(e)\perp_{M} (Id-P)(e)$. Let $v$ and $w$ denote the tripotents given by Proposition~\ref{prop1}$(c)$. Since $M$ and $N$ are non-trivial, $e\in \partial_e \left( \mathcal{B}_W \right)$ and $W = M\oplus^{\ell_{\infty}} N$ it follows that $\|P(e)\| = \|(Id-P)(e)\| =1$. \smallskip
	
\textbf{Claim 1}: $e = v + w$. \smallskip

Indeed, if $v =0$, we have $F_e^{N} =  N_{\circ} \cap  \{ e \}_{_{'}}= \{ v \}_{_{'}} = \emptyset$, and hence $\{ w \}_{_{'}} = F_e^{M} =  M_{\circ} \cap  \{ e \}_{_{'}} =  \{ e \}_{_{'}}$. Theorem~\ref{t ER nclosed faces predual ball} proves that $e = w$. Similarly, if $w=0$ we get $e = v$.  We can therefore assume that $v,w\neq 0$, and hence the norm-closed faces $\{ v \}_{_{'}}$ and $\{ w \}_{_{'}}$ are non-trivial.\smallskip 

On the one hand, since $\{ v \}_{_{'}} = F_e^{N} \subseteq \{ e\}_{_{'}}$ and $\{ w \}_{_{'}} = F_e^{M} \subseteq \{ e\}_{_{'}}$, it follows that $v \leq e$ and $w \leq e$ (cf. Theorem~\ref{t ER nclosed faces predual ball}), and thus $v+w \leq e$.\smallskip 

On the other hand, thanks to Proposition~\ref{prop1}(b), $\{ e\}_{_{'}} = \hbox{conv}(F_e^{N} \cup F_e^{M}) = \hbox{conv}(\{ v\}_{_{'}} \cup \{ w\}_{_{'}})$. Given $\psi \in \{ e\}_{_{'}}$, there exist a real number $t \in [0,1],$ $\phi_{1} \in F_e^{N}=\{ v \}_{_{'}}$ and $\phi_{2} \in F_e^{M}= \{ w \}_{_{'}}$ such that $\psi = t\phi_{1} +(1-t) \phi_{2}$. Hence, $\psi(v+w)= t\phi_{1}(v+w) +(1-t) \phi_{2}(v+w) = t \phi_{1}(v) + (1-t) \phi_{2}(w) = 1,$ where we applied that $v$ and $w$ are orthogonal. The arbitrariness of $\psi\in \{e\}_{_{'}}$ implies that $\{ e\}_{_{'}} \subseteq \{ v+w\}_{_{'}},$ and thus $e \leq v+w$. Therefore $e=v+w$.\smallskip

\textbf{Claim 2}: $P(e) =v$ and $(Id-P)(e) = w$.\smallskip 

For each $\phi \in \{ w \}_{_{'}} = \{ e\}_{_{'}} \cap M_{\circ} $, we have $\phi \mid_{M} \equiv 0$ and $\| \phi \| = 1 = \phi(w) = \phi(e)$. It then follows that $\phi ((Id-P)(e)) = \phi(e) = 1$ as $P(e) \in M$. This proves that $\{ w \}_{_{'}} \subseteq \{ (Id-P)(e) \}_{_{'}}$. It follows from Lemma~\ref{lemma4} that 
$$ (Id-P)(e) = w + x_{0}  $$
with $x_{0} = P_{0}(w)((Id-P)(e)) \in W_{0} (w)$.\smallskip

Similarly, from $\{ v \}_{_{'}} = \{ e \}_{_{'}} \cap N_{\circ}$, we get $\{ v \}_{_{'}} \subseteq \{ P(e) \}_{_{'}} $ and Lemma~\ref{lemma4} proves that
$$  P(e) = v + y_{0} $$   
with $y_{0} = P_{0}(v)(P(e)) \in  W_{0} (v)$. Combining the previous two identities we conclude that
$$ e=v+w = P(e) + (Id-P)(e) = v+ y_{0} +w +x_{0}.$$

Since $e$ is complete, $W_{0}(e) = \{0\}$ and $W =  W_{2}(e) \oplus W_{1}(e)$. From Claim 1 we know that $v+w = e =v+ y_{0} +w +x_{0}$, that is $x_{0}+y_{0} = 0$ with $x_{0}\in W_{0} (w)$ and $y_{0} \in W_{0} (v)$. Since $v\perp w$, the Peirce projections associated with $v$ and $w$ commute, that is, $P_j (v) P_k (w) =P_k (w) P_j (v)$ for all $j,k\in \{0,1,2\}$. Then, by applying $P_{0}(w)$ at both sides, and noting that $P_{0}(w)(y_{0}) = P_0(w) P_0(v) P(e) =   P_0(v) P_0(w) P(e) \in W_{0} (w) \cap W_{0} (v) \subseteq W_{0} (v+w) = \{0\}$, we deduce that 
$$ P_{0}(w)(x_{0}) + P_{0}(w)(y_{0}) = x_{0} + 0 = 0. $$
We have therefore shown that $x_{0} = y_{0} = 0,$ and thus $P(e) = v$ and $(Id-P)(e) = w $ are two orthogonal  tripotents in $W$.\smallskip

\textbf{Claim 3}: $P(e) =v \in \partial_{e}(\mathcal{B}_{M})$ and $(Id-P)(e) =w \in \partial_{e}(\mathcal{B}_N)$.\smallskip

We shall only prove that $P(e) \in \partial_{e}(\mathcal{B}_{M})$, the other statement follows by analogous arguments. Suppose, on the contrary, that there are two different elements $x,y \in \mathcal{B}_{M}$ and $t \in (0,1)$ such that $P(e) = t x + (1-t) y$. Clearly, $x + (Id-P)(e) = x + w$ and $y + (Id-P)(e) = y + w$ both lie in $\mathcal{B}_{W}$ because $(Id-P)(e) = w$ has norm-one and $W = M\oplus^{\ell_{\infty}} N$. The identity
$$t (x + (Id-P)(e)) + (1-t) (y + (Id-P)(e)) = e ,$$ contradicts that $e$ is a complete tripotent in $W$. 
\end{proof}

The previous proposition admits the following interesting corollary which proves an algebraic property of the extreme points of the closed unit balls of any non-trivial $M$-summands in a real JBW$^*$-triple. 

\begin{corollary} \label{coro1}
Let $W$ be a real JBW$^*$-triple. Suppose additionally that $M$ and $N$ are two non-trivial $M$-summands of $W$ with $W= M \oplus^{\ell_\infty} N$. Suppose we take $v \in \partial_{e}(\mathcal{B}_{M})$ and $w \in \partial_{e}(\mathcal{B}_N)$, then $v$ and $w$ are two orthogonal tripotents in $W$ and $v \pm w$ is a complete tripotent in $W$. 
\end{corollary}

\begin{proof} It is well known that $v + w \in \partial_{e}(\mathcal{B}_{W}) = \partial_{e}(\mathcal{B}_{M}) \mathop{\oplus}\limits^{\ell_{\infty}}  \partial_{e}(\mathcal{B}_{N})$. It immediately follows from Theorem~\ref{t complete tripotents and faces} that $P(v+w) = v$ and $(Id-P)(v+w) = w$ are two orthogonal tripotents in $W$.  
\end{proof}

The next observation will also play a key role in our arguments. 

\begin{corollary}\label{c one summand is trivial}  Let $W$ be a real JBW$^*$-triple. Suppose additionally that $M$ and $N$ are two $M$-summands of $W$ with $W= M \oplus^{\ell_\infty} N$, and let $P$ denote the natural projection of $W$ onto $M$. Let $e\in W$ be a complete tripotent {\rm(}equivalently, each extreme point of the closed unit ball of $W${\rm)}. Then $P(e) =0$ implies that $M=\{0\}$ and $(Id-P) (e) =0$ implies $N=\{0\}$. 
\end{corollary}

\begin{proof} If $M,N\neq \{0\},$ Theorem~\ref{t complete tripotents and faces} implies that $\|P(e)\| = \|(Id-P) (e)\|=1$, which contradicts the corresponding assumption. 
\end{proof}

\section{Characterization of $M$-ideals in real JB$^*$-algebras}\label{kernel real JBWstar algebras}

The setting of real JB$^*$-algebras offers some more favourable conditions to understand $M$-ideals, and it is worth to be considered by itself.  Recall that every real JB$^*$-algebra is a real JB$^*$-triple with respect to the triple product defined in \eqref{eq triple product JBstar}. In particular all results in the previous section apply in this setting. More concretely, let $A$ be a real JBW$^*$-algebra with unit denoted by $\mathbf{1}$, and denote by $A_{skew}$ and $A_{sa}$, the set of all skew and self-adjoint elements in $A$, respectively. Suppose that $M$ and $N$ are two $M$-summands of $A$ with $A = M\oplus^{\ell_{\infty}} N$ (note that $M$ and $N$ must be weak$^*$-closed, see Remark~\ref{r M-summands are weak*-closed}). Let $P$ stand for the natural projection of $A$ onto $M$.  Since $\mathbf{1}$ is a complete tripotent in $A$, by Proposition~\ref{prop1} and Theorem~\ref{t complete tripotents and faces}, $P(\textbf{1}) $ and $(Id-P) (\textbf{1}) $ are two orthogonal tripotents in $A$ and $$\{\mathbf{1} \}_{_{'}} = \hbox{conv} \left(\{ P(\textbf{1}) \}_{_{'}} \cup \{ (Id-P) (\textbf{1}) \}_{_{'}} \right).$$ The particular meaning of the partial ordering among tripotents in $A$ implies that $p=P(\textbf{1}) $ and $q=(Id-P) (\textbf{1}) $ are in fact projections in $A$ with $p+q =\mathbf{1}$. \smallskip

Recall that the set of all normal states of $A$ is precisely the set $$\{\mathbf{1} \}_{_{'}} = \{\varphi\in A_* : \varphi \geq 0, \ \|\varphi\|= 1\}  = \{\varphi\in A_* : \varphi (\mathbf{1})= \|\varphi\|= 1\},$$ which coincides with the set of all normal states of the JBW-algebra $A_{sa}$. Contrary to the complex case, the set $\{\mathbf{1} \}_{_{'}} $ does not separate the points the points of $A$, however it is a norming set for $A_{sa}$, because the latter is a JBW-algebra (cf. \cite[Corollary 2.17]{AlfShu2003}).\smallskip

According to the usual notation in Jordan algebras, for every $x,a,b \in A$, we denote by $U_{a,b}$ the linear mapping on $A$ given by $U_{a,b} (x) = \{ a,x^*,b  \}$, we simply write $U_{a}$ for $U_{a,a}$.

\begin{proposition}\label{prop3}  Let $A$ be a real JBW$^*$-algebra with unit $\mathbf{1}$. Suppose that $M$ and $N$ are two $M$-summands of $A$ with $A = M\oplus^{\ell_{\infty}} N$. Let $P$ be the projection from $A$ onto $M$, and let $p$ and $q$ denote the orthogonal projections $P(\textbf{1})$ and $(Id-P) (\textbf{1}),$  respectively. Then the following statements hold: \begin{enumerate}[$(a)$] 
		\item $k \circ M \subseteq M$ and $k \circ N \subseteq N$ for every $k \in A_{skew}$.
		\item $p \circ \frac{m + m^*}{2}= \frac{m + m^*}{2}$ and $q \circ \frac{m + m^*}{2} = 0$ for every $m \in M$.
		\item $q \circ \frac{n + n^*}{2} = \frac{n + n^*}{2}$ and $p \circ \frac{n + n^*}{2} = 0$ for every $n \in N$. 
	\end{enumerate}
\end{proposition}

\begin{proof} $(a)$ Note that by Corollary~\ref{c image under an inner derivation}, for any $a,b \in A$, we have $\delta(a,b)(M) \subseteq M$ and $\delta(a,b)(N) \subseteq N,$ $\delta(a,b) = L(a,b) - L(b,a)$. Consequently, for any $k \in A_{skew}$  (i.e., $k = -k^{*}$) and $m\in M$ we get  
$$M\ni  \delta (k,\mathbf{1})(m) = \{ k,\mathbf{1},m \} - \{\mathbf{1}, k, m \} = 2(k \circ\mathbf{1}) \circ m = 2 k \circ m ,$$ which implies that $k \circ m = \delta(k,\mathbf{1}) \left(\frac{1}{2} m\right) \in M$. Similarly, $k \circ n \in N$ for every $n \in N$.\smallskip

$(b)$ We claim that for $x \in A$, the condition $\phi (x) = 0$ for every $\phi \in \{ q \}_{_{'}}$ implies $p \circ (x + x^*) = (x + x^*).$\smallskip 

Fix an arbitrary  $\psi \in \{ p \}_{_{'}}.$ Note that for every $z \in A$, $U_{p}(z) = \{ p,z^*,p \} = 2 (p \circ z) \circ p - p  \circ z $. Since $\psi \in \{ p \}_{_{'}},$ $\psi = \psi P^{1} (p) = \psi P_{2} (p) = \psi L(p,p) = \psi L(p,\mathbf{1}) = \psi Q(\mathbf{1}) $ (cf. \cite[Lemma 2.7]{PeSta2001}), we have  $$\begin{aligned}
\psi (z)&= \psi (p \circ z ) = \psi (U_{p}(z)) = \psi \left(U_{p}\left(\frac{z+z^*}{2}\right)\right)\\
&= \psi \left(p \circ \frac{z+z^*}{2} \right)= \psi \left(\frac{z+z^*}{2} \right),
\end{aligned}$$ for all $z\in A$. This implies, in particular, that \begin{equation}\label{eq psi vanishes}\begin{aligned}
		\psi(q \circ (z+ z^*))  &= \psi(U_{p} (q \circ (z+ z^*))) \\
		&= \psi \{ p,\{ q,q, z+z^* \},p  \}  =0, \ \hbox{ for all } z\in A,  \psi \in \{ p \}_{_{'}},
	\end{aligned}
\end{equation} where to deduce the last equality we applied that the (self-adjoint) element $\{ q,q, z+z^* \}$ lies in $A_{2}(q) \oplus A_{1}(q),$ and since $p\perp q$, we derive from Peirce arithmetic that $\{ p,\{ q,q, z+z^* \},p  \} \in A_{0-2+0} (q) \oplus A_{0-1+0} (q) =\{0\}$. \smallskip

Take now $x \in A$ satisfying $\phi (x) = 0$ for every $\phi \in \{ q \}_{_{'}}$. For any $\varphi \in \{\mathbf{1} \}_{_{'}} = \hbox{conv} (\{ p \}_{_{'}} \cup \{ q \}_{_{'}})$, there exist a real number $t \in [0,1]$, $\phi \in  \{ q \}_{_{'}}$ and $\psi \in \{ p \}_{_{'}}$ such that $ \varphi = t \phi + (1-t) \psi $ and hence $$ \varphi (q \circ (x+x^*)) = t \phi (q \circ (x+x^*)) + (1-t ) \psi (q \circ (x+x^*)) =0 ,$$
since, by assumptions, $\phi(q \circ (x+x^*)) = \phi (x+x^*) = 0$, and by \eqref{eq psi vanishes}, $\psi (q \circ (x+x^*))= 0$. Therefore, $\varphi (  (\mathbf{1} -p) \circ (x+x^*)  ) = \varphi (q \circ (x+x^*)) = 0$ for all $\varphi \in \{\mathbf{1} \}_{_{'}}$. Since the normal states in $A$ (i.e. functionals in $\{\mathbf{1}\})_{_{'}}$) separate the points in $A_{sa}$, we get $(\mathbf{1} -p) \circ (x+x^*) = 0,$ equivalently $p \circ (x+x^*) = x+x^*$. This finishes the proof of the claim. \smallskip

Take now $m\in M$. By hypotheses, $\phi (m) =0$ for all $\phi\in \{q\}_{_{'}}$. Therefore the claims gives the statement in $(b)$.\smallskip

$(c)$ The proof of this statement is similar to the one given for $(b)$.
\end{proof}

For each element $a$ in a Jordan algebra $\mathcal{J}$, we shall denote by $M_a$ the (Jordan) multiplication operator by $a$ defined by $$ M_a (x) = a\circ x, \ \ (x\in \mathcal{J}).$$ Elements $a$ and $b$ in $\mathcal{J}$ are said to \emph{operator commute} in $\mathcal{J}$ if the multiplication operators $M_a$ and $M_b$ commute in $B(\mathcal{J})$, i.e., $$(a\circ x) \circ b = a\circ (x\circ b), \hbox{ for all $x$ in $\mathcal{J}$.}$$ The centre of $\mathcal{J}$, $Z(\mathcal{J})$, is the set of all elements $z$ in $\mathcal{J}$ such that $z$ and $b$ operator commute for every $b$ in $\mathcal{J}$. \smallskip 

We can now establish our main conclusion concerning $M$-summands of real JBW$^*$-algebras. 

\begin{theorem}\label{t M.-ideals in real JBWstar algebras}  Let $A$ be a real JBW$^*$-algebra, and let $M$ be a closed subspace of $A$. Then $M$ is an $M$-summand in $A$ if and only if it is a weak$^*$-closed {\rm(}Jordan{\rm)} ideal of $A$. 
\end{theorem}

\begin{proof} It is known that every weak$^*$-closed Jordan ideal of a real JBW$^*$-algebra is an $M$-summand (see Remark~\ref{r wstar triple ideals are M-summands} and the comments prior to it).\smallskip 
	
To prove the necessary condition, let $A$ be a real JBW$^*$-algebra with unit $\mathbf{1},$ and let $M$ be an $M$-summand in $A$ which is the image of an $M$-projection $P$ on $A$. Setting $N:=(Id-P) (A)$, it follows that $M$ and $N$ are weak$^*$-closed (see Remark~\ref{r M-summands are weak*-closed}) and $A = M\oplus^{\ell_{\infty}} N$. Let us denote $p=P(\textbf{1})\in M$ and $q=(Id-P) (\textbf{1})\in N.$ We know that $p$ and $q$ are two orthogonal projections in $A$ satisfying $p+q = \mathbf{1}$, $\{\mathbf{1} \}_{_{'}} = \hbox{conv} (\{ p \}_{_{'}} \cup \{ q \}_{_{'}})$ (see Proposition~\ref{prop1} and Theorem~\ref{t complete tripotents and faces}).\smallskip %We also know from Proposition~\ref{prop3} that the following statements hold:
% \begin{enumerate}[$(a)$] 
 %	\item $k \circ M \subseteq M$ and $k \circ N \subseteq N$ for every $k \in A_{skew}$.
 %	\item $p \circ \frac{m + m^*}{2}= \frac{m + m^*}{2}$ and $q \circ \frac{m + m^*}{2} = 0$ for every $m \in M$.
 %	\item $q \circ \frac{n + n^*}{2} = \frac{n + n^*}{2}$ and $p \circ \frac{n + n^*}{2} = 0$ for every $n \in N$. 
%	\end{enumerate}
	
\noindent We shall first prove that 
\begin{equation}\label{eq p and q are central} \hbox{  $p$ and $q$ are central projections in $A$}.
\end{equation}
	
To simplify the notation, fix arbitrary $m \in M$, $n \in N$ and denote $m= \frac{1}{2}(m + m^*) + \frac{1}{2} (m - m^*) = h + k$, $n =\frac{1}{2}(n + n^*) + \frac{1}{2} (n- n^*)= \tilde{h} + \tilde{k},$ where $h, \tilde{h} \in A_{sa}$ and $k, \tilde{k} \in A_{skew}$. Observe that \begin{equation}\label{eq p circ k} p \circ k =(p \circ k)\circ \mathbf{1} = (p \circ k) \circ q + (p \circ k) \circ p, 
\end{equation} and the decomposition is unique. Clearly  $p \circ k, k\in A_{skew}$, and hence it follows from Proposition~\ref{prop3}$(a)$ that $p \circ k \in M$, $(p \circ k) \circ q \in N$ and $(p \circ k) \circ p \in M$. it follows from the identity in \eqref{eq p circ k} that $(p \circ k) \circ q = 0$. We can similarly obtain $(q \circ k) \circ p = 0$. By combining this conclusion with Proposition~\ref{prop3}$(b)$ and $(c)$ we get \begin{equation}\label{eq p and q are central for M} \begin{aligned}
	(p \circ m) \circ q &= (p \circ (h+k)) \circ q  = (p \circ h) \circ q +  (p \circ k) \circ q\\
	&= 0 + k\circ q =0, \\ 
	(q \circ m) \circ p &= (q \circ h) \circ p +(q \circ k) \circ p = 0, \hbox{ for all } m\in M.
\end{aligned}
\end{equation} If in the above arguments we replace $m\in M$ with $n\in N$ we obtain 
\begin{equation}\label{eq p and q are central for N} (p \circ n) \circ q = (q \circ n) \circ p = 0, \hbox{ for all } n\in N.
\end{equation} Now, by combining \eqref{eq p and q are central for M}, \eqref{eq p and q are central for N} and $A = M\oplus^{\ell_\infty} N$ we deduce that \begin{equation}\label{eq p and q are central for A} (p \circ a) \circ q = (q \circ a) \circ p = 0, \hbox{ equivalently, } U_{p,q} (a) =0, 
\end{equation} for all $a\in A$. This implies that 
$$A = U_{\mathbf{1}}(A) = U_{p+q}(A)= U_{p}(A) \oplus U_{p,q}(A) \oplus U_{q}(A) = U_{p}(A) \mathop{\oplus}\limits^{\ell_{\infty}} U_{q}(A),$$
and consequently $p$ and $q$ are central projections in $A$. It is perhaps worth to note that in the last equality we wrote ``$\mathop{\oplus}\limits^{\ell_{\infty}}$'' because $p\perp q$.\smallskip

Since $p$ and $q$ are central orthogonal projections in $A$ with $$A = U_{p}(A) \mathop{\oplus}\limits^{\ell_{\infty}} U_{q}(A) = A_{2}(p) \mathop{\oplus}\limits^{\ell_{\infty}} A_{2}(q),$$ the predual of $A$ decomposes in the form 
$$ A_{*} = (A_{2}(p))_{*} \mathop{\oplus}\limits^{\ell_{1}}  (A_{2}(q))_{*} ,$$ that is, ${P_{2}(p)}_{*}$ and ${P_{2}(q)}_{*}$ (i.e. the transposed maps of the corresponding Pierce 2-projections associated with $p$ and $q$ restricted to $A_{*}$) are $L$-projections on $A_*$. Since $L$-projections on a real or complex Banach space commute (cf. \cite[Lemma 2.2]{Cunningham1960}), we have $[P_{*}, {P_{2}(p)}_{*}] = 0 = [P_{*}, {P_{2}(q)}_{*}]$ where $P_{*}$ is the transposed of $P$ restricted to $A_*$, and hence \begin{equation}\label{eq P anad P2p P2q commute}  [P, P_{2}(p)] = 0= [P_{2}(q), P].
\end{equation}

In this case, since $p$ and $q$ are central we also have $P_2 (p) (a) = p\circ a$ and $P_2 (q) (a) = q\circ a$, for all $a\in A$. So, it follows from \eqref{eq P anad P2p P2q commute} that $$ M\ni P(p\circ a) =  P P_2(p) (a) =P_2(p) P(a) = p\circ P(a) = p\circ m,$$  for all $a = m+n \in A = M\oplus N.$ If we we combine this conclusion with Proposition~\ref{prop3}$(a)$ and $(b)$ we obtain $$M\ni p\circ m = p\circ \frac{m +m^*}{2} + p\circ \frac{m -m^*}{2}, \hbox{ and }  p\circ \frac{m -m^*}{2}\in M,$$ and thus  $ \frac{m +m^*}{2}= p\circ \frac{m +m^*}{2}\in M$ for all $m\in M$. By observing that (by Proposition~\ref{prop3}$(a)$ and $(b)$) $$\begin{aligned}
M\ni m &= (p+q) \circ m \\
  &= p\circ \frac{m +m^*}{2} + p\circ \frac{m -m^*}{2} +  q\circ \frac{m +m^*}{2} + q\circ \frac{m -m^*}{2} \\
 &= p\circ \frac{m +m^*}{2} + p\circ \frac{m -m^*}{2} +  0 + q\circ \frac{m -m^*}{2}
\end{aligned}$$ with $ p\circ \frac{m +m^*}{2} + p\circ \frac{m -m^*}{2}\in M$ and $q\circ \frac{m -m^*}{2}\in N$, we conclude that $q\circ \frac{m -m^*}{2} =0$, and thus $q \circ m = 0$ for all $m\in M$. We similarly get $p \circ N = \{0\}.$ Both equalities show that $$ M = p\circ M = P_2 (p) (M) \subseteq P_2(p) (A) = p\circ A = p \circ M + p \circ N = p \circ M,$$ and $$ N = q\circ N = P_2 (q) (M) \subseteq P_2(q) (A) = q\circ A = q \circ M + q \circ N = p \circ N,$$ that is $M = p\circ A = P_2(p) (A)$ and $N= q\circ A = P_2 (q) (A)$ are two (Jordan) ideals of $A$.\end{proof}

The description of all $M$-ideals of a real JB$^*$-algebra can be now obtained as a corollary.

\begin{corollary}\label{c M-ideals of real JB-algebras}
Let $M$ be a subspace of a real JB$^*$-algebra $A$. Then $M$ is an $M$-ideal of $A$ if and only if it is a {\rm(}Jordan{\rm)} ideal of $A$.
\end{corollary}

\begin{proof} Suppose $M$ is an $M$-ideal of $A$. In this case, there exists an $L$-projection $Q$ on $A^*$ whose image is $M^{\circ}$, and $A^{*} = M^{\circ} \oplus^{\ell_1} (Id-P) (A^*)$. It is known that $Q^*$ is an $M$-projection on $A^{**}$ and $A^{**} = Q^* (A^{**})  \oplus^{\ell_\infty} (Id-Q)^* (A^{**}) = \overline{M}^{w^*} \oplus^{\ell_\infty} (Id-Q)^* (A^{**})$, where $\overline{M}^{w^*}$ denotes the weak$^*$-closure of $M$ in $A^{**}$. Theorem~\ref{t M.-ideals in real JBWstar algebras} asserts that $\overline{M}^{w^*}$ is a weak$^*$-closed ideal of $A^{**}$, and in particular, $M$ is an ideal of $A$.\smallskip
	
Finally, every Jordan ideal of $A$ is an $M$-ideal by Remark~\ref{r wstar triple ideals are M-summands} and the comments prior to it.
\end{proof}

The description of all $M$-ideals in an arbitrary real C$^*$-algebra is just a straightforward consequence of the previous Corollary~\ref{c M-ideals of real JB-algebras} by just observing that every Jordan ideal of a real C$^*$-algebra is a two-sided ideal.   

\begin{corollary}\label{c M-ideals of real Cstar-algebras}
	Let $M$ be a subspace of a real C$^*$-algebra $A$. Then $M$ is an $M$-ideal of $A$ if and only if it is a two-sided ideal of $A$.
\end{corollary}

\section{Main results: Characterization of $M$-ideals in real JB$^*$-triples}\label{kernel real JBWstar triples}

The goal of this section is to study $M$-summands of real JBW$^*$-triples, and as a consequence to characterize the $M$-ideals of a real JB$^*$-triple.\smallskip

In this section $W$ will be a real JBW$^*$-triple, $M$ will be an $M$-summand of $W$, that is, $W = M\oplus^{\ell_{\infty}} N$, where $N$ is another closed subspace of $W$ {\rm(}both $M$ and $N$ are automatically weak$^*$-closed, see Remark~\ref{r M-summands are weak*-closed}{\rm)}. We can clearly assume that $M$ and $N$ are both non-trivial. The symbol $P$ will stand for the $M$-projection on $W$ whose image is $M$, and we shall fix a complete tripotent $e$ in $W$. By Theorem~\ref{t complete tripotents and faces} (and Proposition~\ref{prop1}) the elements $v=P(e)$ and $w=(Id-P)(e)$ are two orthogonal non-zero tripotents in $W$ with $e = v +w$, $N_{\circ} \cap  \{ e \}_{_{'}}= \{ v \}_{_{'}} \perp_{L} M_{\circ} \cap \{ e \}_{_{'}} =\{ w \}_{_{'}}$, and $\hbox{conv}(\{ v \}_{_{'}}\cup \{ w \}_{_{'}}) = \{ e \}_{_{'}}$. Clearly, $v$ and $w$ are projections in the real JBW$^*$-algebra $W_2(e)$ (and in the JBW-algebra $W^{1} (e)$).\smallskip
 
For every $m \in M^{**}$ with Peirce decomposition $m = P_1(e) (m) + P_2(e)(m)$ and $P_2(e) (m) = P^{1} (e) (m) + P^{-1} (e) (m)$ in which $P^{1} (e) (m) \in \left(W_{2}(e)\right)_{sa}$ and $k \in \left(W_{2}(e)\right)_{skew}$. A similarly for elements in $N$. We keep this notation in the rest of this section. \smallskip

Our first goal is an appropriate extension of Proposition~\ref{prop3} and Theorem~\ref{t M.-ideals in real JBWstar algebras}.

\begin{lemma}\label{l L(e,k)} Let $e$ be a tripotent in a real JB$^*$-triple $E$. Then for each $k$ in $E^{-1} (e)$ {\rm(}i.e. $k^{*_e} = \{e,k,e\} = -k${\rm)} we have $L(e,k) = - L(k,e)$ and hence $\delta(e,k) = 2L(e,k) =-2 L(k,e)$. 
\end{lemma}

\begin{proof} Consider any $ x \in E$ with $x = P_0 (e) (x) + P_{1} (e) (x) + P_{2} (e) (x)$. By linearity, it suffices to prove that $L(e,k) (P_j (e) (x)) = - L(k,e) (P_j (e) (x))$ for all $j\in \{0,1,2\}$. Since $e,k\in E_2(e)\perp E_0(e)$, it follows that $L(e,k)(P_0(e)(x)) = 0 = - L(k,e)(P_0(e)(x))$. \smallskip
	
Having in mind that $e,k,P_2(e)(x)\in E_2 (e)$, and the triple product on $E_2(e)$ is uniquely determined by the Jordan product $\circ_e$ and the involution $*_{e}$ of the unital real JB$^*$-algebra $E_2(e)$, we get 
$$\begin{aligned}
L(e,k)(P_2(e)(x)) &= \{e,k,P_2(e)(x)\}= (e\circ_e k^{*_e}) \circ_e P_2(e)(x)  \\
&+ (P_2(e)(x) \circ_e k^{*_e}) \circ_e e - (e\circ_e P_2(e)(x))\circ_{e} k^{*_e} \\
&= - (e\circ_e k) \circ_e P_2(e)(x)  - (P_2(e)(x) \circ_e k) \circ_e e \\
&+ (e\circ_e P_2(e)(x))\circ_{e} k \\
&= - (e^{*_e}\circ_e k) \circ_e P_2(e)(x)  - (P_2(e)(x) \circ_e k) \circ_e e^{*_e} \\
&+ (e^{*_e}\circ_e P_2(e)(x))\circ_{e} k \\
&=-L(k,e) (P_2(e) (x)).
\end{aligned}$$
 	 
We consider next the Peirce-$1$ subspace. By the Jordan identity and Peirce arithmetic, we have
	\begin{equation*}
		\begin{split}
			\frac12 L(e,k) (P_1(e)(x))&= L(e,k) \{ e,e, P_1(e)(x)\} \\
			&= \{ L(e,k) (e), e, P_1(e)(x) \} - \{ e,L(k,e)(e), P_1(e)(x) \} \\
			&+\{ e,e, L(e,k) (P_1(e)(x)) \} \\
			&=  \{ L(e,k)(e),e, P_1(e)(x) \} - \{ e,k, P_1(e)(x) \} \\
			&+ \frac{1}{2}  L(e,k) (P_1(e)(x)).  
		\end{split}
	\end{equation*}
Since $ L(e,k)(e) = k^{*_e} = -k = -  L(k,e)(e)$, the above identity shows that 
	$$  L(e,k) (P_1(e)(x)) =  -\{ k,e, P_1(e)(x) \} = -L(k,e) (P_1(e)(x)), $$
	which concludes the proof.
\end{proof}

\begin{proposition}\label{prop4} Under the assumptions stated at the beginning of this section we have:  
\begin{enumerate}[$(a)$]
\item $L(e,k)(M) = -L(k,e)(M) \subseteq M$ and $L(e,k)(N) = -L(k,e)(N) \subseteq N,$ for every $k\in W^{-1} (e)$. 
\item  $\{P^{1} (e) (m), e, v \} = P^{1} (e) (m)$ and $\{ P^{1} (e) (m), e, w \} = 0,$ for every element $m\in M$.
\item  $\{P^{1} (e) (n), e, w \} = P^{1} (e) (n)$ and $\{ P^{1} (e) (n), e, v \} = 0,$ for every $n\in N$.
\item $v$ and $w$ are central projections in $W_{2}(e)$, equivalently, $$W_1(v) \cap W_1 (w)=\{0\}, \hbox{ or  } W_2 (e) = W_2(v)\oplus^{\ell_{\infty}} W_2 (w).$$
\end{enumerate} 
\end{proposition}

\begin{proof} $(a)$ By combining Lemma~\ref{l L(e,k)}, Corollary~\ref{c image under an inner derivation} and the fact that $M$ and $N$ are $M$-summands of $W$, we get $$\begin{aligned}
2 L(e,k)(M) &=-2 L(k,e)(M) = \delta(e,k) (M) \subseteq M, \\
2 L(e,k)(N) &= -2 L(k,e)(N) = \delta(e,k) (M) \subseteq N
\end{aligned}$$

$(b)$ Inspired by the arguments we employed in the proof of Proposition~\ref{prop3}(b), for any $\psi \in \{ w \}_{_{'}}$, we define a real linear functional $\psi_{w}\in W_*$ given by $\psi_{w} (x) := \psi \{ w,e,x \}$ for any $x \in W$. By observing that 
$$1=\| \psi(w) \|  \leq \| \psi_{w} \| \leq \| \psi \| = 1, $$
we have $\psi_{w} \in \{ w \}_{_{'}} = \{ e \}_{_{'}} \cap M_{\circ}$ and $\psi_{w} (m) = 0$ since $m \in M$. Then, $$\psi \{ w,e, P_{1}(e) (m) \} +\psi \{ w,e, P_{2}(e) (m) \}  = \psi \{ w,e, m \} = \psi_{w}(m)= 0$$
for all $\psi \in \{ w \}_{_{'}}.$ It follows from $\psi \in \{ w \}_{_{'}}\subseteq \{ e \}_{_{'}}$, that $\psi  = \psi P_{2}(e) = \psi P^{1}(e)$ (cf. \cite[Lemma 2.7]{PeSta2001}). Since, by Peirce arithmetic, $\{ w,e, P_{1}(e) (m) \} \in W_{1}(e)$, we conclude from the previous identity that $\psi \{ w,e, P_{1}(e) (m) \} = 0$. Therefore,  
\begin{equation}\label{eq psi P1 e m is zero} \begin{aligned}
		0= \psi_{w}(m) &= \psi_{w}(P_{2}(e)(m)) = \psi_{w}(P^{1}(e)(m)) \\
		&= \psi(P^{1}(e) (m)) \ \ \ (\forall m\in M, \psi\in \{ w \}_{_{'}}),
	\end{aligned}
\end{equation}
where the third we applied that $\psi_{w}\in \{ e \}_{_{'}}$ and hence $\psi_{w} = \psi_w P^1 (e)= \psi_w P_2 (e)$ (see \cite[Lemma 2.7]{PeSta2001}). We compute next the values of $\psi$ at an element of the form $\{P^{1}(e)(m),v, v \}$. Concretely, having in mind that $\psi\in \{ w \}_{_{'}}$ we deduce from \eqref{eq positive normal states and hermitian values} that 
\begin{equation}\label{eq psi at P1em v v} \psi \{P^{1}(e)(m),v, v \} = \psi\{w,  \{P^{1}(e)(m),v, v \},w\} = 0,
\end{equation}	because, by Peirce arithmetic and the fact that that $w\perp v$ and hence $w\in W_0 (v),$ we derive that 
$$  \{w,  \{P^{1}(e)(m),v, v \},w\} \in \{W_0(v), W_{1}(v)\oplus W_2 (v), W_0(v)\} =  \{0\}.$$

Take now $\phi \in  \{ v \}_{_{'}}\subseteq  \{ e \}_{_{'}}$. Since $\phi = \phi P^1 (e) =\phi P_2(e)$ and $e = v+ w$ with $w\perp v$, we derive that \begin{equation}\label{eq phi P1em e v zero}
	\begin{aligned}
		\phi \{ P^{1}(e)(m),e, v \} &= \phi \{ P^{1}(e)(m), e, e \} \\
		&= \phi \left(P_2(e) +\frac12 P_1(e) \right) P^{1}(e)(m ) 
		\\ &= \phi P^{1}(e)(m) = \phi (m).
	\end{aligned}
\end{equation}

Let us recall that $ \{ e \}_{_{'}} = conv(\{ v \}_{_{'}} \cup \{ w \}_{_{'}})$ (see Theorem~\ref{t complete tripotents and faces}). So, for any $\varphi \in \{ e \}_{_{'}} = \hbox{conv} (\{ v \}_{_{'}} \cup \{ w \}_{_{'}})$, there exist a real number $t \in [0,1]$, $\phi \in  \{ v \}_{_{'}}$ and $\psi \in \{ w \}_{_{'}}$ such that 
$$ \varphi = t \phi + (1-t) \psi .$$
It then follows from \eqref{eq phi P1em e v zero}, \eqref{eq psi P1 e m is zero}, and   \eqref{eq psi at P1em v v} that 
\begin{equation*}
\begin{split}
  \varphi(\{ P^{1}(e)(m), e, v \} - P^{1}(e)(m)) &= t \phi(\{ P^{1}(e)(m),e, v \} - P^{1}(e)(m)) \\
  &+(1-t) \psi(\{ P^{1}(e)(m),e, v \} - P^{1}(e)(m) \}) \\
      &=  (1-t) \psi \{ P^{1}(e)(m),e, v \} \\
      &=  (1-t) \psi \{P^{1}(e)(m),v, v \}= 0.  
\end{split}
\end{equation*}

By the arbitrariness of $\varphi \in \{ e \}_{_{'}}$ in the previous identity combined with the fact that $\{ P^{1}(e) (m),e, v \} - P^{1}(e) (m)$ is a self-adjoint elements in $W_{2} (e)$, we conclude that $\{ P^{1}(e)(m),e, v \} - P^{1}(e) (m) = 0$. Therefore, $\{ P^{1} (e)(m), e, v \} = P^{1} (e) (m),$ and $\{ P^{1}(e)(m) ,e, w \} = 0$ as desired.\smallskip

The proof of $(c)$ follows by similar arguments.\smallskip

$(d)$  Having in mind that $P^{-1}(e) (m), \{ P^{-1} (e)(m),e,v \}\in W^{-1} (e)$, the conclusion in statements $(a)$ and $(b)$ give 
$$\begin{aligned}
M\ni \{ P^{-1} (e)(m),e,v \} &=  \{ \{ P^{-1} (e)(m),e,v \}, e, e \} \\ &= \stackrel{\in M}{\overbrace{\{ \{ P^{-1} (e)(m),e,v \}, e, v \}}} +   \stackrel{\in N}{\overbrace{\{ \{ P^{-1} (e)(m),e,v \}, e, w \}}}
\end{aligned},$$ which implies that \begin{equation}\label{eq 0212} 0 =\{ \{ v,e,P^{-1} (e)(m) \}, e, w \} = \left( P^{-1} (e)(m)\circ_e v\right)\circ_e w.
\end{equation} We can similarly obtain that \begin{equation}\label{eq 0212 b} \left( P^{-1} (e)(m)\circ_e w\right)\circ_e v =\{ \{ w,e,P^{-1}(e) (m) \}, e, v \} =0.
\end{equation} Similar identities hold when $m\in M$ is replaced by any $n\in N$.\smallskip 

The desired conclusion in $(d)$ will follow if we prove that $\{ v, W , w \} = \{ v, W_2(e) , w \} = 0$ (the first equality follows from Peirce arithmetic). Since $W = M\oplus^{\ell_{\infty}} N$, it suffices to show that $\{ v, M, w \} = \{ v, N, w \}  = 0$. For $m\in M$ we have  
\begin{equation}\label{vmv}
\{ v,m,w \} = \{v, P_{1}(e) (m), w \} + \{ v,P^{1}(e) (m), w \} + \{ v,P^{-1}(e) (m) ,w \},
\end{equation}
with $\{v, P_{1} (e) (m), w \} = 0 $ by Peirce arithmetic. Note that $v, w ,P^{1}(e) (m) $ and $P^{-1}(e) (m)$ all lie in the real JBW$^*$-algebra $W_{2} (e)$, and hence by $(b)$, $v\perp w$, \eqref{eq 0212} and \eqref{eq 0212 b} we arrive at 
\begin{equation*}
  \begin{split}
     &\{ v, P^{1}(e) (m), w \} + \{ v, P^{-1}(e) (m) ,w \} = (v \circ_e P^{1}(e) (m)) \circ_e w \\
     & + (w \circ_e P^{1}(e) (m)) \circ_e v - (v \circ_e w) \circ P^{1}(e) (m) - (v \circ_e P^{-1}(e) (m)) \circ_e w \\
       & - (w \circ_e P^{-1}(e) (m)) \circ_e v + (v \circ_e w) \circ_e P^{-1}(e) (m) =0. % \\
       %& = P^{1}(e) (m) \circ_{e} w =\{P^{1} (e), e, w\} =  0.
  \end{split}
\end{equation*}
Therefore, back to \eqref{vmv}, we conclude that $\{ v,m,w \} = 0$ for every $m \in M$. The proof of $\{v, N, w\} =0$ follows by very close arguments. \smallskip

We have therefore shown that $\{v,W,w\}=\{0\},$ which in particular implies that $Q(v,w) =0.$ Recall that $w$ and $v$ are two orthogonal projections in $W_2(e)$ with $v+ w =e$ (and hence two compatible tripotents in $W$) and the well known decomposition $$W_2 (e) = W_2(v) \oplus W_2 (w) \oplus \Big(W_1(v) \cap W_1 (w)\Big),$$ with $W_2 (v) \subseteq W_0(w)$ and $W_2 (w) \subseteq W_0(v)$. Since $$\begin{aligned}
P_2 (e) &= Q(e)^2 = Q(v+w)^2 = (Q(v)+Q(w)+ 2 Q(v,w) )^2 = (Q(v)+Q(w))^2 \\
&= Q(v)^2 + Q(w)^2 + Q(v) Q(w) +  Q(w) Q(v)\\ 
&= P_2(v) + P_2(w),  
\end{aligned}$$ we obtain $W_1(v) \cap W_1 (w) =\{0\}$ and $W_2(e) = W_2 (v) \oplus^{\ell_{\infty}} W_2(w).$
\end{proof}

In our last technical step we shall prove that $v$ and $w$ actually induce an $M$-decomposition of the whole $W$.

\begin{proposition}\label{prop5} Under the assumptions stated at the beginning of this section the subspaces $W_{1,0} = W_{1} (v) \cap W_{0} (w)$ and $W_{0,1} = W_{0} (v) \cap W_{1} (w)$ are orthogonal in $W$, and  
$$W = \left(W_{1} (v) \oplus W_{2} (v)\right) \mathop{\oplus}\limits^{\ell_{\infty}}  \left( W_{1} (w) \oplus W_{2} (w)\right).$$
Observe that $(P_{1}(v) +P_{2}(v))$ and $(P_{1}(w) +P_{2}(w))$ are two $M$-projections on $W$ whose images are the $M$-summands $W_{1} (v) \oplus W_{2} (v)$ and $W_{1} (w) \oplus W_{2} (w)$, respectively. Furthermore, if $P$ denotes the $M$-projection of $W$ onto the $M$-summand $M$, the projections $P$, $(Id-P),$ $(P_{1}(v) +P_{2}(v))$ and $(P_{1}(w) +P_{2}(w))$ are pairwise commuting.
\end{proposition}

\begin{proof} Let us begin by proving the first statement. By Proposition~\ref{prop4}$(d)$, $W_2 (e) = W_2(v)\oplus^{\ell_{\infty}} W_2 (w)$ (in particular $W_1(v) \cap W_1 (w) =\{0\}$). Fix $a \in W_{1,0}$ and $b \in W_{0,1}$. By Peirce arithmetic \begin{equation}\label{eq 3a 312} \{a,a,b\}\in W_{1-1+0} (v) \cap W_{0-0+1} (w) = W_0 (v)\cap W_{1} (w) = W_{0,1},\end{equation} and
\begin{equation}\label{eq 3b 312}  \{w,b,a\} \in W_{0-0+1} (v)\cap W_{2-1+0} (w) = W_1 (v)\cap W_1 (w)=\{0\}.
\end{equation}
 By assumptions $b= 2 \{ w,w, b \}$, and hence by applying that $a \perp w,$ \eqref{eq 3b 312} and Peirce arithmetic we arrive at
\begin{equation*}
  \begin{split}
     \{ a,a,b \} &= 2 \{ a,a, \{ w,w, b \} \} = 2 \{ a,a, L(b,w)(w)\}  \\
       &= 2 (L(b,w)\{ a,a,b \}-\{ L(b,w)a,a,b \} +\{ a,L(w,b)a,b \} ) \\
       &= 2 (L(b,w)\{ a,a,b \} - 0 + 0 ) \\
       &= 2 L(b,w)\{ a,a,b \},
  \end{split}
\end{equation*} where by \eqref{eq 3a 312} $$\begin{aligned}
L(b,w)\{ a,a,b \}&\in W_{0-0+0} (v) \cap  W_{1-2+1} (w) = W_0(v) \cap W_0 (w)\\
&\subseteq W_{0} (v+w) =W_0 (e)=\{0\}.
\end{aligned}$$ Therefore, $\{a,a,b\}=0,$ or equivalently $a \perp b$ for every $a\in W_{1,0}$ and $b \in W_{0,1},$ that is, $W_{1,0} \perp W_{0,1}$. \smallskip

Since under our assumptions and what we have just proved we get $$W_1 (e) = \left( W_0(v) \cap W_1(w) \right) \oplus^{\ell_{\infty}} \left( W_1(v) \cap W_0(w) \right) = W_{0,1} \oplus^{\ell_{\infty}} W_{1,0} $$  (recall that orthogonality implies $M$-orthogonality \cite[Lemma 1.3$(a)$]{FriedmanRusso1985}).\smallskip
 
Since $v$ and $w$ are orthogonal, and thus compatible, the Peirce projections associated with $v$ and $w$ commute, which implies that
$$\begin{aligned}
W_1 (v) &=  \Big(W_1 (v)\cap W_2(w)\Big) \oplus \Big(W_1 (v)\cap W_1(w)\Big) \oplus \Big(W_1 (v)\cap W_0(w)\Big) \\
&\subseteq  \Big(W_1 (v)\cap W_0 (v)\Big) \oplus \{0\} \oplus \Big(W_1 (v)\cap W_0(w)\Big) =W_{1,0} \subseteq W_1 (v) \\ 	
\end{aligned},$$ and similarly
$W_1 (w)  = W_{0,1}.$ The last two identities and the first statement assure that $W_2 (v)\oplus W_1 (v) = W_2 (v)\oplus W_{1,0} \perp W_2(w)\oplus W_{0,1}= W_2(w)\oplus W_1(w),$ and thus $W_2 (v)\oplus W_1 (v)$ and  $W_2(w)\oplus W_1(w)$ are $M$-orthogonal. \smallskip

Summarizing all the previous conclusions we deduce that $$\begin{aligned}
W &= W_2(e) \oplus W_1 (e) = W_2(v) \oplus W_2(w) \oplus W_1(v) \oplus W_1 (w) \\
&= (W_{1} (v) \oplus W_{2} (v)) \mathop{\bigoplus}\limits^{\perp, \ell_{\infty}}  (W_{1} (w) \oplus W_{2} (w)).
\end{aligned}$$ Observe that $W(v) = W_{1} (v) \oplus W_{2} (v)$ and $W(w) = W_{1} (w) \oplus W_{2} (w)$ are two orthogonal weak$^*$-closed triple ideals of $W$, whose direct sum is the whole $W$, the $M$-projections of $W$ onto $W(v)$ and $W(w)$ are $(P_{1}(v) +P_{2}(v))$ and $(P_{1}(w) +P_{2}(w))$, respectively. The corresponding predual spaces, $W(v)_{*}$ and $W(w)_{*},$ are $L$-summands of $W_*$ with associated $L$-projections $(P_{1}(v) +P_{2}(v))_*$ and $(P_{1}(w) +P_{2}(w))_*$, respectively. We apply once again the mentioned result by Cunningham (cf. \cite[Lemma 2.2]{Cunningham1960}) to conclude that $P_*$ commutes with $(P_{1}(v) +P_{2}(v))_*$ and $(P_{1}(w) +P_{2}(w))_*$, and hence $P$ commutes with $(P_{1}(v) +P_{2}(v))$ and $(P_{1}(w) +P_{2}(w))$. 
\end{proof}

Let us recall that the set of all \emph{contractive perturbations} of a subset $\mathcal{S},$ of the closed unit ball, $\mathcal{B}_{X}$, of a Banach space $X$ is defined by $$\hbox{\rm cp}(\mathcal{S})= \hbox{\rm cp}_{X} (\mathcal{S}) := \{ x\in X : \|x \pm s \| \leq 1, \hbox{ for all } s\in \mathcal{S} \}\subseteq \mathcal{B}_{X}.$$ For each $a\in \mathcal{B}_{X}$ we write cp$(a)$ for the set cp$(\{a\})$. The geometric characterization of tripotents in real or complex JB$^*$-triples established in \cite[Theorem 2.3]{FerMarPe2004geometric} asserts that a norm one element $e$ in a real JB$^*$-triple $E$ is a tripotent if, and
only if, the sets $$ D_{1} (e) := \{ y \in {E} : \hbox{ there	exists } \alpha>0 \hbox{ with } \| e \pm \alpha y \| = 1 \}$$ and
$$ D^{\prime}_{2} (e) := \{ y \in {E} : \| x + \beta y \| = \max \{1,\|\beta y\|\} \hbox{ for all } \beta\in \mathbb{R} \}$$ coincide. It is actually shown in the proofs of \cite[Theorems 2.1 and 2.3]{FerMarPe2004geometric} that for each tripotent $e$ in $E$ we have $D_1 (e) =E_0 (e)$, which in terms of the set of all contractive perturbations of $e$ it can be restated as follows:   
\begin{equation}\label{eq contractive perturbation of a tripotent}\hbox{\rm cp}_{E}(\{e\})\subseteq D_1(e)
	\cap E_1 = \mathcal{B}_{E_0 (e)} = \{e\}^{\perp}\cap \mathcal{B}_{E} 
	\subseteq \hbox{\rm cp}_{E} (\{e\}).
\end{equation} This was already observed in the case of complex JB$^*$-triples, with almost identical arguments, in \cite[identity $(6)$ in page 360]{FerMarPe2012}. \smallskip

We are now in a position to state the main result of this paper. 

\begin{theorem}\label{t M-summands in real JBW*-triples are ideals} Let $M$ be a closed subspace of a real JBW$^*$-triple $W$. Then $M$ is an $M$-summand of $W$ if and only if it is a weak$^*$-closed triple ideal of $W$.
\end{theorem}

\begin{proof} We have seen in Remark~\ref{r wstar triple ideals are M-summands} that every weak$^*$-triple ideal is an $M$-summand.\smallskip

Suppose that $M$ is an $M$-summand of $W$. Let $P$ denote the $M$-projection of $W$ onto $M$, and let $N= (Id-P)(W)$. Then $M$ is weak$^*$-closed (cf. Remark~\ref{r M-summands are weak*-closed}).  We can clearly assume that $M$ is non-trivial. By Theorem~\ref{t complete tripotents and faces} (and Proposition~\ref{prop1}) the elements $v=P(e)$ and $w=(Id-P)(e)$ are two orthogonal non-zero tripotents in $W$ with $e = v +w$. If we chain the results in Propositions~\ref{prop4} and \ref{prop5}, we deduce that the real JBW$^*$-subtriples  $W(v) = W_{1} (v) \oplus W_{2} (v)$ and $W(w) = W_{1} (w) \oplus W_{2} (w)$ are two orthogonal weak$^*$-closed triple ideals of $W$, whose direct sum is the whole $W$. As before, by Cunningham's result \cite[Lemma 2.2]{Cunningham1960}, the projections $P,$ $Id-P,$ $(P_{1}(v) +P_{2}(v))$ and $(P_{1}(w) +P_{2}(w))$ pairwise commute.\smallskip

Let us consider the real JBW$^*$-triple $W(v)$ and the subspaces $M_1 = W(v)\cap M$ and $N_1 = W(v) \cap N$. Since the projections $P,$ $Id-P,$ $(P_{1}(v) +P_{2}(v))$ and $(P_{1}(w) +P_{2}(w))$ pairwise commute, it can be easily deduced that $W(v) = M_1 \oplus^{\ell_{\infty}} N_1,$ that is, $M_1$ and $N_1$ are two $M$-summands of $W(v)$. Furthermore, the corresponding $M$-projections of $W(v)$ onto $M_1$ and $N_1$ are $\tilde{P} =P (P_{1}(v) +P_{2}(v))$ and $Id-\tilde{P} = (Id-P) (P_{1}(v) +P_{2}(v))$, respectively. Having in mind that $v$ is a complete tripotent in $W(v)$ with $\tilde{P} (v) = P (P_{1}(v) +P_{2}(v)) (v) = P(v) =v$ and $(Id-\tilde{P}) (v) = (Id-P) (P_{1}(v) +P_{2}(v)) (v) =0$, Corollary~\ref{c one summand is trivial} implies that $N_1=\{0\}.$ Consequently, $N \subseteq W(w)$ and $W(v) = M_1 = W(v)\cap M\subseteq M$.  We can similarly obtain that $M\cap W(w) = \{0\}$, and thus $M \subseteq W(v)$ and $N \subseteq W(w)$. Therefore $N = W(w)$ and $M = W(v)$, which gives the desired result.\smallskip

Alternatively, if there exists a non-zero element $y\in N_1 = W(v) \cap N\subseteq W(u)$, which can assumed to have norm-one, by observing that $v$ is a complete tripotent in $W(v)$ (equivalently, $\left(W(v)\right)_{0}(v)=\{0\}$), $v = P(e)\in M$, and $W = M\oplus^{\ell_\infty} N$, we get $\| y \pm v \| = max \{ \|y\|, \|v \|  \} = 1$. By \eqref{eq contractive perturbation of a tripotent} (cf. \cite[identity $(6)$ in page 360]{FerMarPe2012}), we obtain that $y \in \hbox{cp}_{_{W(v)}} ( v) = \left(W(v)\right)_{0}(v) =\{0\}$, contradicting that $y$ is non-zero. Similarly, $W(w)\cap M =\{0\}.$
\end{proof}

As a consequence of our previous theorem we can now characterize $M$-ideals of real JB$^*$-triples as triple ideals.

\begin{theorem}\label{thm1}
Let $E$ be a real JB$^*$-triple. Then the $M$-ideals of $E$ are precisely the {\rm(}triple{\rm)} ideals of $E$.
\end{theorem}

\begin{proof} We only need to prove that every $M$-ideal is triple ideal (cf. Remark~\ref{r wstar triple ideals are M-summands}). Suppose $M$ is an $M$-ideal of $E$. We consider $M$ as a closed subspace of the real JBW$^*$-triple $E^{**}$. It is known that under this assumptions, $M^{**}\cong \overline{M}^{w^*}  \cong (M_{\circ})^{\circ}$ is an $M$-summand of $E^{**}$. Theorem~\ref{t M-summands in real JBW*-triples are ideals} implies that $M^{**}\cong \overline{M}^{w^*}$ is a weak$^*$-closed triple ideal of $E^{**},$ and hence $M$ is a triple ideal of $E$.
\end{proof}

\begin{corollary}\label{c complexificaiton of M-ideals} Let $M$ be a closed subspace of a real JB$^*$-triple $E$, and let $\mathcal{E}$ denote the JB$^*$-triple obtained by complexifying $E$. Then $M$ is an $M$-ideal of $E$ if and only if the natural complexification of $M$ is an $M$-ideal in $\mathcal{E}$.
\end{corollary}

Let us recall that a \emph{real TRO} is a closed linear subspace $Z \subseteq B(K, H),$ for real Hilbert spaces $K$ and $H$, satisfying $Z Z^* Z \subseteq Z$ (see, for example, \cite{BleLeMerdy2004, Sharma14}). Clearly, every real TRO is a real JB$^*$-triple, and thus the following corollary is an immediate consequence of Theorem~\ref{thm1}.

\begin{corollary}\label{c real TRO}
Let $Z$ be a real TRO. Then the $M$-ideals of $Z$ coincide with its {\rm(}triple{\rm)} ideals.
\end{corollary}

\noindent\textbf{Acknowledgements}\smallskip

\noindent Third author supported by Junta de Andalucía grant FQM375, MCIN (Ministerio de Ciencia e Innovación, Spain)/ AEI/10.13039/501100011033 and “ERDF A way of making Europe”  grant PID2021-122126NB-C31, MOST (Ministry of Science and Technology of China) grant G2023125007L, and by the IMAG--Mar{\'i}a de Maeztu grant CEX2020-001105-M/AEI/10.13039/ 501100011033. Fourth author supported by China Scholarship Council.

%---------------------------------

\subsection*{Data availability}

There is no data associate for the submission entitled ``$M$-ideals, yet again: the case of real JB$^*$-triples''.

\subsection*{Statements and Declarations}

The authors declare they have no financial nor conflicts of interests.

%\reference here
%\bibliographystyle{siam}
%\bibliographystyle{unsrt}
%\bibliographystyle{abbrv}
%\bibliography{sn-a}

% ------------------------------------------------------------------------
\end{document}